\documentclass[12pt]{article}
\usepackage{amssymb}
\usepackage{dsfont}
\usepackage{calrsfs} 
\usepackage{mathtools}
\usepackage{hyperref}
\usepackage{xcolor}
\hypersetup{
	colorlinks,
	linkcolor={magenta!100!black},
	citecolor={blue!100!black},
	urlcolor={purple!100!black}
}
\newtheorem{theorem}{Theorem}
\newtheorem{proposition}{Proposition}[theorem]
\newtheorem{remark}{Remark}
\newtheorem{coro}{Corollary}
\newtheorem{lemma}{Lemma}
\newcommand{\E}{\mathbb{E}}
\begin{document}
%
%
%\centerline{\huge Nonparametric local linear regression estimation under censorship model}
\noindent{\huge Strong consistency of the  nonparametric local linear regression estimation under censorship model}
\vskip 3mm

\vskip 5mm
\noindent {\bf F. Bouhadjera}

\noindent Universit\'e Badji-Mokhtar, Lab. de Probabilit\'es et Statistique. BP 12, Annaba, 23000, Alg\'erie. 

\noindent \url{ferielbouhadjera@gmail.com}

\noindent {\bf E. Ould Sa\"id}

\noindent Universit\'e du  Littoral C\^ote d'Opale, Lab. de Math. Pures et Appliqu\'ees, IUT de Calais. 19, rue Louis David. Calais, 62228, France

\noindent \url{elias.ould-said@univ-littoral.fr}

\noindent {\bf M. R. Remita}

\noindent Universit\'e Badji-Mokhtar, Lab. de Probabilit\'es et Statistique. BP 12, Annaba, 23000, Alg\'erie. 

\noindent \url{medriad.remita@univ-annaba.dz}

\vskip 3mm

\noindent {\bf Abstract:} We introduce and study a local linear nonparametric regression estimator for censorship model. The main goal  of this paper is, to establish the uniform almost sure consistency result with rate over a compact set for the new estimate. To support our theoretical result, a simulation study has been done to make comparison with the classical regression estimator.

\noindent {\bf Key Words:} censored data; local linear estimation; rate of consistency; regression function; survival analysis; uniform almost sure consistency.
\vskip 3mm 
\vskip 4mm

\section{Introduction}

Nonparametric regression is a smoothing method for recovering a regression function from data, which has no restriction on its form. In this context, the modelization of the relationship between two random variables has been studied by many researchers in complete and incomplete data. The most used approaches are kernel methods (see e.g: \hyperref[Guessoum2008]{Guessoum and Ould Sa\"id (2008)}). It is well known that the local linear smoothing procedure has many desirable advantages, for an extensive discussion for what regards the bias and boundary effects in theses issues, we can refer to  \hyperref[fan1992]{Fan (1992)}, \hyperref[fan1996]{Fan and Gijbels (1996)} for univariate case and \hyperref[fan2003]{Fan and Yao (2003)} for the multivariate case.\\
In this work, we deal with the problem of the estimation of the regression function under right censorship by using the local linear approach. In this framework, we define a local linear regression estimator by taking in account the synthetic data and establish the uniform almost sure consistency with rate of the resulting estimator. Among the studies dedicated to local linear fit in the case of censored data and without pretending to exhaustivity, we quote \hyperref[cai2003]{Cai (2003)} who proposed an estimator of the regression function based on the generalization of the weighted least squares method. We point out that his study does penalize the censoring survival function. Recall that \hyperref[fan1994]{Fan and Gijbels (1994)} established the consistency in probability of analogous estimator.  Finally, inspired by the works of \hyperref[beran1981]{Beran (1981)} and \hyperref[dabrowska1987]{Dabrowska (1987)}, \hyperref[kim1998]{Kim (1998)} used the conditional hazard estimator to define a linear local estimator.  We point out that the only result in this framework is the consistency in probability and as far as we know our result is new.\\
The organization of the article is as follows.  In the next section, we give a brief description of nonparametric regression estimates and the proposed methodology. In \hyperref[sect 3]{Section 3}, we provide the hypotheses, main results and a sketch of the proof. A comparative study with the classical kernel estimator with different sizes and censoring rates have been done in \hyperref[sect 4]{Section 4}.  Finally, the proofs are postponed to the \hyperref[sect 5]{Section 5}.
%----------------------- Model ----------------------------------%
\section{Model and estimators}\label{sect 2}
Consider $n$ independent and identically distributed replications of a couple $(X_i,Z_i)$ having  the same  distribution as the pair $(X,Z)$ where $Z$ is the  interest random variable (r.v.) with unknown distribution function  $H$ and $X$ is the corresponding covariate. Under this setting, the purpose of this paper is to consider a regression model $Z=\mu(X)+\varepsilon$ with
\begin{equation*}
\mu(x)=\mathbb{E}[Z|X=x]=\frac{\displaystyle  \int_{\mathbb{R}} z \;f(z, x)\;dz}{f(x)}=:\frac{S_0(x)}{f(x)}
\end{equation*}
where $f(\cdot,\cdot)$, $f(\cdot)$ are the joint density of $(X,Z)$ and marginal density of $X$ respectively and the white noise  $\varepsilon$ is assumed to be normally distributed with zero mean and variance $\sigma_\varepsilon^2$. In the situation of censored data, we do not observe $Z$ but only $Y=\min(Z,C)$ and $\delta=\mathds{1}_{\{Z\leq C\}}$ where $C$ is the censoring variable. In what follows we assume that $Z$ and $C$ are independent. This assumption is required to ensure the identifiability of the model.
Under this setting, we consider a specific transformation of the data that take into account the effect of the censoring in the distribution: the so-called synthetic data introduced by \hyperref[carbonez1995]{Carbonez {\it et al.} (1995)} and used by  \hyperref[kohler2002]{Kohler {\it et al.}(2002)}, \hyperref[Guessoum2008]{Guessoum and Ould Sa\"id (2008)} and a large number of authors given, for $1 \leq i \leq n$, by 
\begin{equation}\label{syntdata}
Z^\star_i=\frac{\delta_i Y_i}{\overline{G}(Y_i)}.
\end{equation}
where $\overline{G}(\cdot)$ denotes the survival function of $C$. In what follows, it is assumed that:
\begin{equation}\label{indep}
(Z_i,X_i)_i \;\; \text{and} \;\; (C_i)_i \;\;\text{are independent.}
\end{equation}
Using the conditional expectation properties and the condition (\ref{indep}), then for all fixed  $x$, we have
\begin{equation*}\label{synt}
\begin{aligned}
\mathbb{E}[Z^\star_1|X_1=x]&=\mathbb{E}\left[\frac{\delta_1 Y_1}{\overline{G}(Y_1)}\big|X_1=x\right]\\
&=\mathbb{E}\left[\mathbb{E}\left[\frac{\mathds{1}_{\{Z_1 \leq C_1\}} Z_1}{\overline{G}(Z_1)}\big|Z_1\right]\big|X_1=x\right]\\
&=\mathbb{E}\left[\frac{Z_1}{\overline{G}(Z_1)}\mathbb{E}\left[\mathds{1}_{\{Z_1 \leq C_1\}}|Z_1\right]X_1=x\right]\\
&=\mathbb{E}[Z_1|X_1=x].
\end{aligned}
\end{equation*}
Now, assume that the second derivative of $\mu(x)$ exists. Based on the approximation of  $\mu(X) \approx \mu(x)+\mu'(x)(X-x)\equiv \alpha+\beta(X-x)$ in a neighborhood of a point $x$, we extend the LLR estimator to the censoring case by substituting $Z$ by $Z^\star$.
The problem of estimating $\mu(x)$ becomes minimizing
\begin{equation}\label{minimiz}
\arg\min_{\alpha,\beta} \sum_{1 \leq i \leq n} \left(Z^\star_i-\alpha-\beta(X_i-x)\right)^2 K_i,
\end{equation}
with $K_i=K\left(h^{-1}(X_i-x)\right) $ for all $i=1,\dots,n$ where $K$ is a kernel  density, and $h:=h_n$ is a  sequence of the the strictly positive numbers which goes to zero as $n$ goes to infinity. By a simple algebra computation, solving (\ref{minimiz}) yields to  the "Pseudo-estimator" defined by 
\begin{equation}\label{LLRE}
\widetilde{\mu}(x)=\frac{\displaystyle\sum_{1 \leq i,j \leq n} w_{i,j}(x) Z^{\star}_j}{\displaystyle \sum_{1 \leq i,j \leq n} w_{i,j}(x)}=:\frac{\widetilde{\mu}_1(x)}{\widehat{\mu}_0(x)}
\end{equation}
where
\[w_{i,j}(x)=(X_i-x)\left((X_i-x)-(X_j-x)\right)K_iK_j.\]
\noindent 	Of course in data analysis, the survival function $\overline{G}(\cdot)$ is unknown and needs to be estimated. This can be done via the \hyperref[kaplan1958]{Kaplan and Meier (1958)} as an estimator of $\overline{G}$ given  by 
\begin{equation}\label{K-M}
\overline{G}_n(t)=
\left\{
\begin{array}{cl}
\displaystyle \prod_{i=1}^{n}{\left(1-\frac{1-\delta_{i}}{n-i+1}\right)}^{\mathds{1}_{\{Y_{i}\leq t\}}} & \quad \text{if} \quad t<Y_{(n)}, \\
0 & \quad \text{otherwise}
\end{array}
\right.
\end{equation}
where $Y_{(1)}\leq Y_{(2)}\leq \dots \leq Y_{(n)}$ are the order statistics of the $Y_{i}$ and $\delta_{i}$ is a concomitant of   $Y_{i}$. The properties of $\overline{G}_n(t)$ have been studied by many authors.
Hence to get a feasible estimator, we replace (\ref{K-M}) in (\ref{syntdata}), then we get 
\begin{equation}\label{syntdata2}
\widehat{Z}^\star_i=\frac{\delta_i Y_i}{\overline{G}_n(Y_i)}  \quad \quad \text{for} \quad 1 \leq i \leq n,
\end{equation}
and substituting (\ref{syntdata2}) in (\ref{LLRE}) we get
\begin{equation}\label{LLR}
\widehat{\mu}(x)=\frac{\displaystyle\sum_{1 \leq i,j \leq n} w_{i,j}(x) \widehat{Z}^{\star}_j}{\displaystyle\sum_{1 \leq i,j \leq n} w_{i,j}(x)}=:\frac{\widehat{\mu}_1(x)}{\widehat{\mu}_0(x)}, \quad \quad \quad \left(\frac{0}{0}=:0\right).
\end{equation}
The estimator $\widehat{\mu}(\cdot)$ is called the local linear regression (LLR) smoother and it has many desirable statistical properties such as avoiding the edge effects (see: ~\hyperref[fan1992]{Fan (1992)}).
%%%%%%%%%%%%%%%%%%%%%%%%%%%%%%%%%%%%%%%%%%%%%%%%%%%%%%%%%%%%%
\begin{remark} If $\beta=0$, we obtain the classical regression (CR) estimator given in \hyperref[Guessoum2008]{Guessoum and Ould Sa\"id (2008)} and 
	defined by 
	\begin{equation}\label{CR}
	\mu_n(x)=\frac{ \displaystyle \sum_{1 \leq j \leq n} \widehat{Z}^\star_j K_j}{  \displaystyle \sum_{1 \leq j \leq n} K_j}
	\end{equation}
\end{remark}
%--------------------- ASSUMPTIONS AND Main results ---------------------------------%
\section{Assumptions and Main results} \label{sect 3}
Let ${\mathcal{C}}_{0}=\left\{x \in \mathbb{R} / f(x) >0 \right\}$  and ${\mathcal{C}}$ be a compact subset of ${\mathcal{C}}_{0}$. Throughout the paper, we assume that for any d.f. $Q$, we have $\tau_Q=\sup\{x,Q(x)<1\}$ the upper endpoint of the support. We assume that $\tau_H>0$ and $0<\overline{G}(\tau_H)<\infty$.\\
When no confusion is possible, we denote by $C$ any generic positive constant. Furthermore, as $Z$ is a lifetime it can be supposed to be bounded. Our assumptions are gathered together for easy references.
\begin{itemize}
	\item[A1.] The bandwidth $h$ satisfies $ \displaystyle\lim_{n \rightarrow \infty} h=0, \quad \lim_{n \rightarrow \infty} n h=+\infty,\quad \lim_{n\rightarrow\infty} \frac{\log n}{nh}=0$.\label{A1}
	\item[A2.] The kernel $K(\cdot)$ is a bounded, symmetric nonnegative function on $\mathcal{C}$. \label{A2}
%	Furthermore, for $j=2,3$ \label{A2}
%	\begin{itemize}
%		\item[i)] $\int t^j K(t)dt<\infty$.\label{A2i}
%		\item[ii)] $\int t^j K^2(t)dt<\infty$.\label{A2ii}
%	\end{itemize}
	\item[A3.] The density function $f(\cdot)$ is continuously differentiable and  $\displaystyle \sup_{x \in \mathcal{C}} |f^{\prime}(x)| <+\infty$. \label{A3}
	\item[A4.] The function $S_0(x)$ is continuously differentiable and $\displaystyle \sup_{x \in \mathcal{C}} |S_0^{\prime}(x)| <+\infty$. \label{A4}
	\item[A5.] The function $\upsilon_k(x)=\int z^k f_{Z,X}(z,x)dz$, %where $f(\cdot,\cdot)$ is the joint density of the couple $(Z,X)$, 
	is continuously differentiable and $\displaystyle \sup_{x \in \mathcal{C}} |\upsilon_k^{\prime}(x)| <+\infty$. \label{A5}
	\item[A6.] There exists $C>0$, $\nu >0$ such that 
	\[\forall (x,y) \in \mathbb{R}^2 \quad \quad \quad |\mu(x)-\mu(y)| \leq C|x-y|^{\nu}.\] \label{A6}
\end{itemize} 

\vspace*{-1cm}
\noindent {\sf Remarks on the Assumptions.}\\
Notices that the assumptions \hyperref[A1]{A1} concerns the bandwidth and is analohous to that is used in   \hyperref[Guessoum2008]{Guessoum and Ould Sa\"id (2008)}. The assumption \hyperref[A2]{A2} deals with the Kernel $K$ and is needed for the convergence of the bias and variance terms. Major standard kernels satisfy these assumptions, for example the Epanechnikov or  Gaussian kernels. The Assumptions \hyperref[A3]{A3},  \hyperref[A4]{A4}, \hyperref[A5]{A5} and \hyperref[A6]{A6} are regularity conditions on  the density $f(\cdot)$, $S_0(\cdot)$, $\mu$ and $\upsilon_k(\cdot)$ respectively.  \vspace*{0.3cm}\\
%%%%%%%%%%%%%%%%%%%%%%%%%%%%%%%%%%%%%%%%%%%%%%%%%%%%%%%%%%%%%%%%%%%%%%%%%%%%%%
The following theorem gives the almost sure (a.s.) consistency of $\widehat{\mu}$ over a compact set $\mathcal{C}$ with rate.
\begin{theorem}\label{theo1}
	Under Assumptions \hyperref[A1]{A1}--\hyperref[A6]{A6}, we have
	\begin{equation*}
	\sup_{x\in {\mathcal{C}}} |\widehat{\mu}(x)-\mu(x)|= \text{O}(h^{\nu})+\text{O}_{a.s.}\left(\sqrt{\frac{\log n}{nh}}\right) \quad \quad \quad \text{as} \quad  n \rightarrow \infty.
	\end{equation*}
\end{theorem}
The proof is based on the following decomposition: 
	\begin{equation*}
	\begin{aligned}
	\widehat{\mu}(x)-\mu(x)&=: \mathcal{B}_1(x)+\frac{1}{\widehat{\mu}_0(x)}\left\{-\mathcal{B}_1(x)\mathcal{B}_2(x) + \mathcal{B}_3(x) +\mathcal{B}_4(x) - \mu(x) \mathcal{B}_2(x)\right\}
	\end{aligned}
	\end{equation*}
	with 
	\begin{eqnarray}
	\mathcal{B}_1(x)&:=& \frac{\mathbb{E}[\widetilde{\mu}_1(x)]}{\mathbb{E}[\widehat{\mu}_0(x)]}-\mu(x), \quad
	\mathcal{B}_2(x):=\widehat{\mu}_0(x)-\mathbb{E}[\widehat{\mu}_0(x)],\nonumber\\ 
	\mathcal{B}_3(x)&:=&\widehat{\mu}_1(x)-\widetilde{\mu}_1(x)\quad \mbox{and}\quad 
	\mathcal{B}_4(x):=\widetilde{\mu}_1(x)-\mathbb{E}[\widetilde{\mu}_1(x)].\nonumber
	\end{eqnarray}
	By triangle inequality, we have
	\begin{equation*}
	\begin{aligned}
	\sup_{x \in \mathcal{C}} |\widehat{\mu}(x)-\mu(x)|& \leq \inf_{x \in \mathcal{C}}  |\mathcal{B}_1(x)| +\frac{1}{\displaystyle\inf_{x \in \mathcal{C}}|\widehat{\mu}_0(x)| }\left\{\sup_{x \in \mathcal{C}}|\mathcal{B}_1(x)\mathcal{B}_2(x)|\right.\\
	&\left.+\sup_{x \in \mathcal{C}}|\mathcal{B}_3(x)|+\sup_{x \in \mathcal{C}}|\mathcal{B}_4(x)|+\sup_{x \in \mathcal{C}}|\mu(x)\mathcal{B}_2(x)|\right\}.
	\end{aligned}
	\end{equation*}
	The proof will be achieved with the following propositions: 
	%%%%%%%%%%%%%%%%%%%%%%%%%%%%%%%%%%%%%
	\begin{proposition}\label{prop1}
		Under Assumptions \hyperref[A1]{A1} and \hyperref[A6]{A6}, we have 
		\begin{equation*}
		\displaystyle \sup_{x \in \mathcal{C}}\left|\mathcal{B}_1(x)\right|=\text{O}\left(h^{\nu}\right) \;\;\;\; \text{as} \;\;\;\; n \rightarrow \infty.
		\end{equation*}
	\end{proposition}
	%%%%%%%%%%%%%%%%%%%%%%%%%%%%%%%%
	\begin{proposition}\label{prop2}
		Under Assumptions \hyperref[A1]{A1-A3}, we have 
		\begin{equation*}
		\sup_{x \in \mathcal{C}}|\mathcal{B}_2(x)|=\text{O}_{a.s.} \left(\sqrt{\frac{\log n }{n h }}\right) \;\;\;\; \text{as} \;\;\;\; n \rightarrow \infty.
		\end{equation*}
	\end{proposition}
    %%%%%%%%%%%%%%%%
   \begin{coro}
	Under the assumption of \hyperref[prop2]{Proposition 2}, there exists a real number $\Gamma>0$ such that: 
	\begin{equation*}
	\sum_{n=1}^{\infty}\mathbb{P} \left( \inf_{x \in \mathcal{C}} \widehat{\mu}_0(x) \leq \Gamma \right) < \infty.
	\end{equation*}
    \end{coro}
	%%%%%%%%%%%%%%%%%%%%%%%%%%%%%%%%
	\begin{proposition}\label{prop3}
		Under Assumptions \hyperref[A1]{A1}, \hyperref[A2]{A2} and \hyperref[A3]{A3}, we have
		\begin{equation*}
		\sup_{x \in \mathcal{C}}\left|\mathcal{B}_3(x)\right|=\text{O}_{a.s.}\left(\sqrt{\frac{\log \log n}{n}}\right) \;\;\;\; \text{as} \;\;\;\; n \rightarrow \infty.
		\end{equation*}
	\end{proposition}
	%%%%%%%%%%%%%%%%%%%%%%%%%%%%%%%%%%%%%%%%
	\begin{proposition}\label{prop4}
		Under Assumptions \hyperref[A1]{A1-A5}, we have
		\begin{equation*}
		\sup_{x \in \mathcal{C}}\left|\mathcal{B}_4(x)\right|=\text{O}_{a.s.}\left(\sqrt{\frac{\log n}{nh}}\right) \;\;\;\; \text{as} \;\;\;\; n \rightarrow \infty.
		\end{equation*}
	\end{proposition}
	%----------------------------SIMULATION PART ------------------------------------------%
\section{Numerical study}\label{sect 4}
\noindent Simulations are conducted to assess the finite sample performance of the proposed estimator $\widehat{\mu}(\cdot)$ given in \hyperref[sect 2]{Section 2} and compare its efficiency and robustness over the classical kernel regression estimator defined in (\ref{CR}). In all the presented curves, we take $K$ to be standard normal density function and the optimal bandwidth $h$ is selected by the well known cross validation method. 
We simulate $n$ points from the following model: $Z_i=X_i+0.2\,\epsilon_i$ where $X_i \leadsto \mathcal{N}(0,1)$ and $\epsilon_i \leadsto \mathcal{N}(0,1)$. The censoring time is distributed as $C_i \leadsto \mathcal{N}(c,1)$ where $c$ is a constant that adjusts the censoring percentage (C.P.). We compute the transformed data obtained via (\ref{synt}) where the K.M. estimator is defined in (\ref{K-M}).
\begin{figure}[!h]
	\begin{minipage}[c]{.26\linewidth}
		\includegraphics[height=2in, width=2.3in]{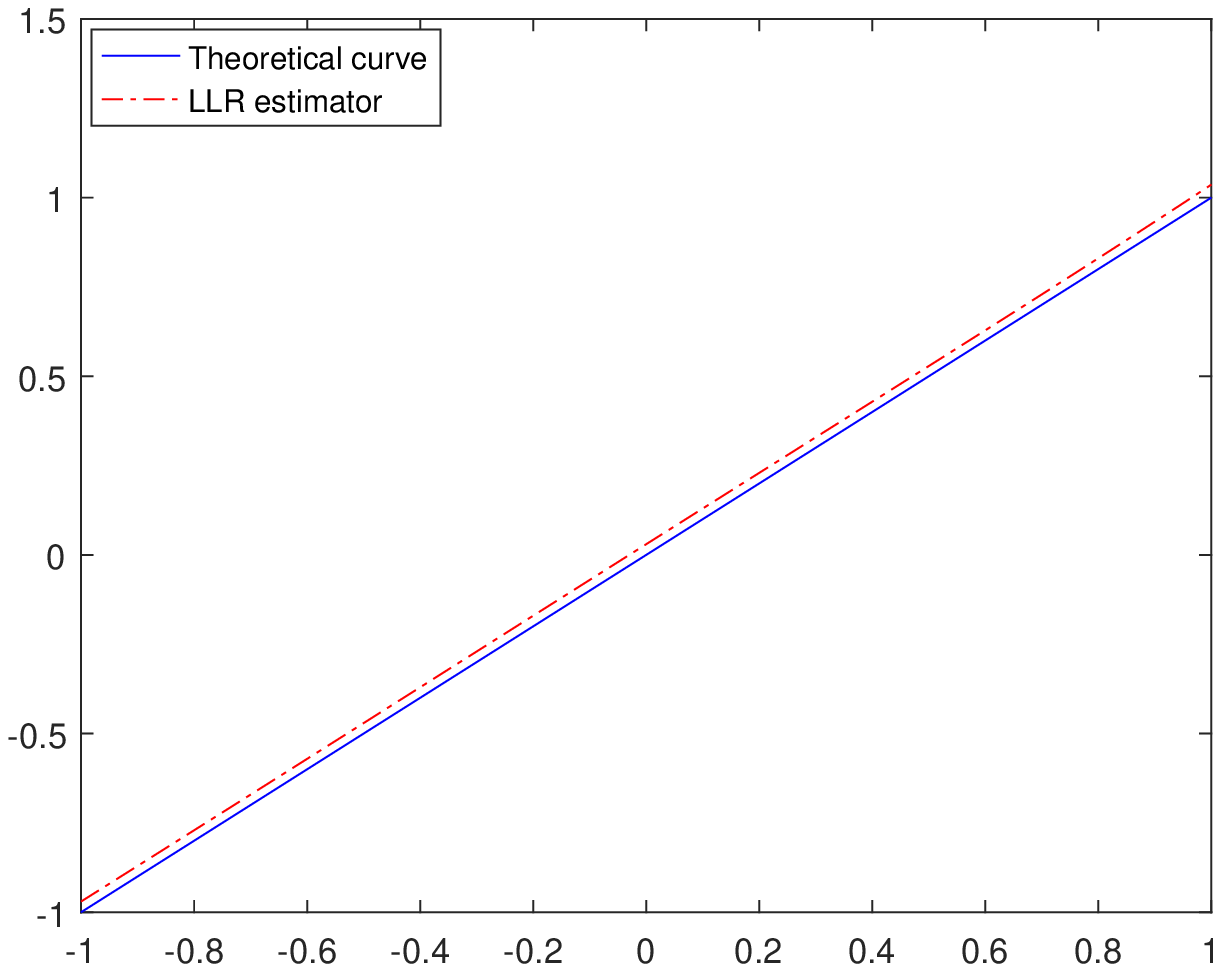}
	\end{minipage} \hfill
	\begin{minipage}[c]{.26\linewidth}
		\includegraphics[height=2in, width=2.3in]{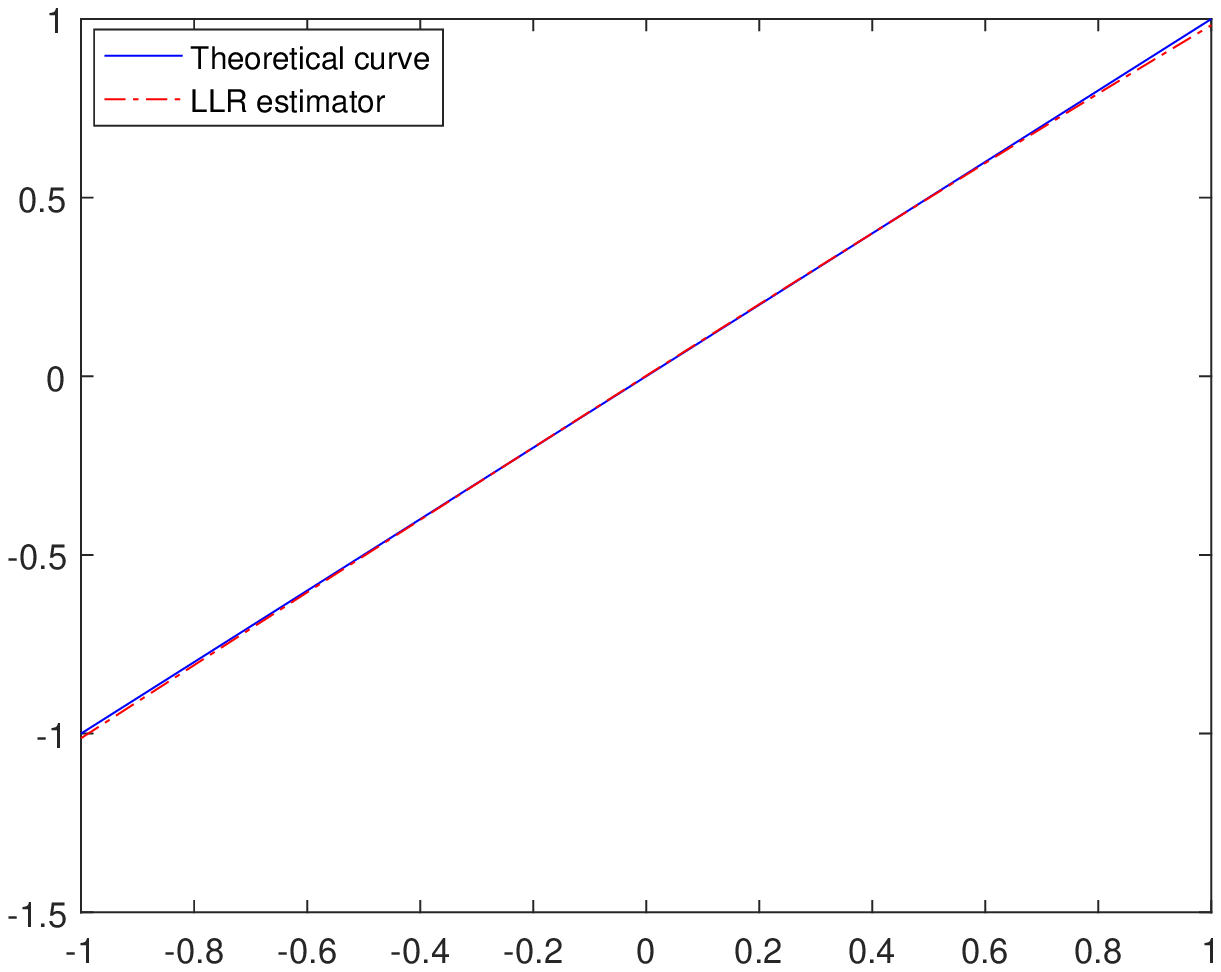}
	\end{minipage} \hfill
	\begin{minipage}[c]{.26\linewidth}
		\includegraphics[height=2in, width=2.3in]{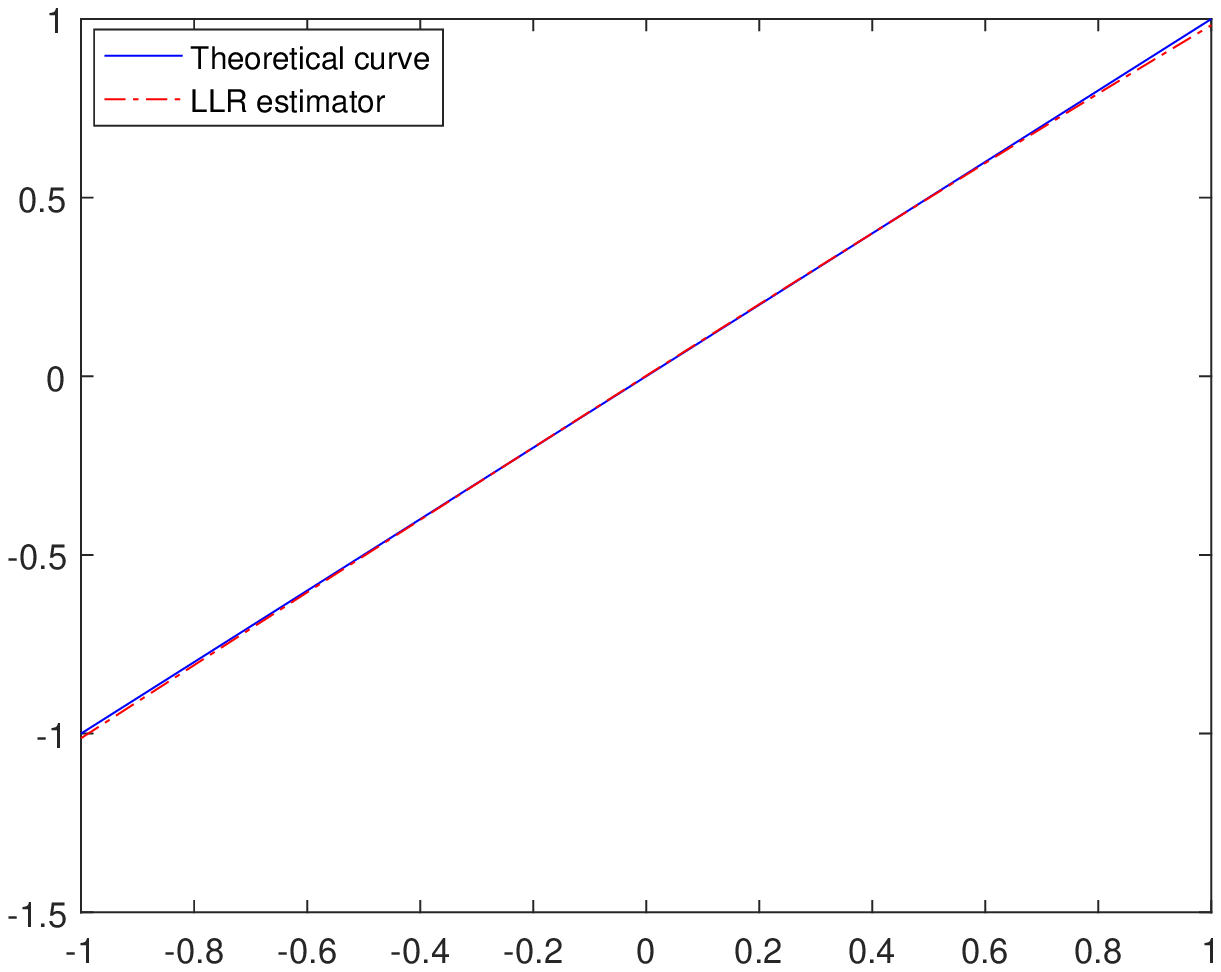}
	\end{minipage}\hfill\hfill
	\caption{\textcolor{blue}{$\quad\mu(\cdot)$}, \textcolor{red}{$\widehat{\mu}(\cdot)$} with C.P.$ \approx30\%$ for $n=100, 300$ and $500$ respectively.}\label{figure1}
\end{figure}
%%%%%%%%%%%%%%%
\begin{figure}[!h]
	\begin{minipage}[c]{.26\linewidth}
		\includegraphics[height=2in, width=2.3in]{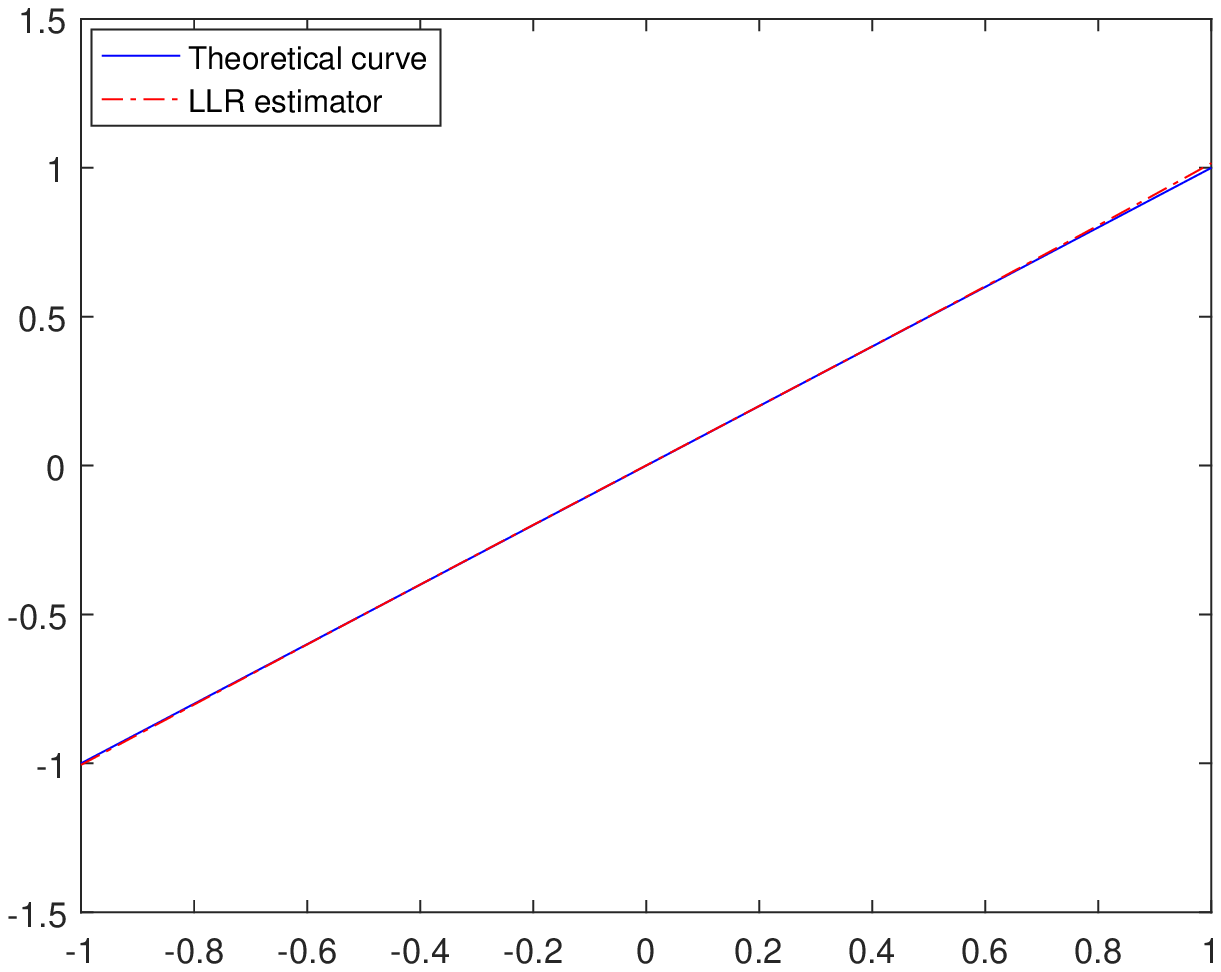}
	\end{minipage} \hfill
	\begin{minipage}[c]{.26\linewidth}
		\includegraphics[height=2in, width=2.3in]{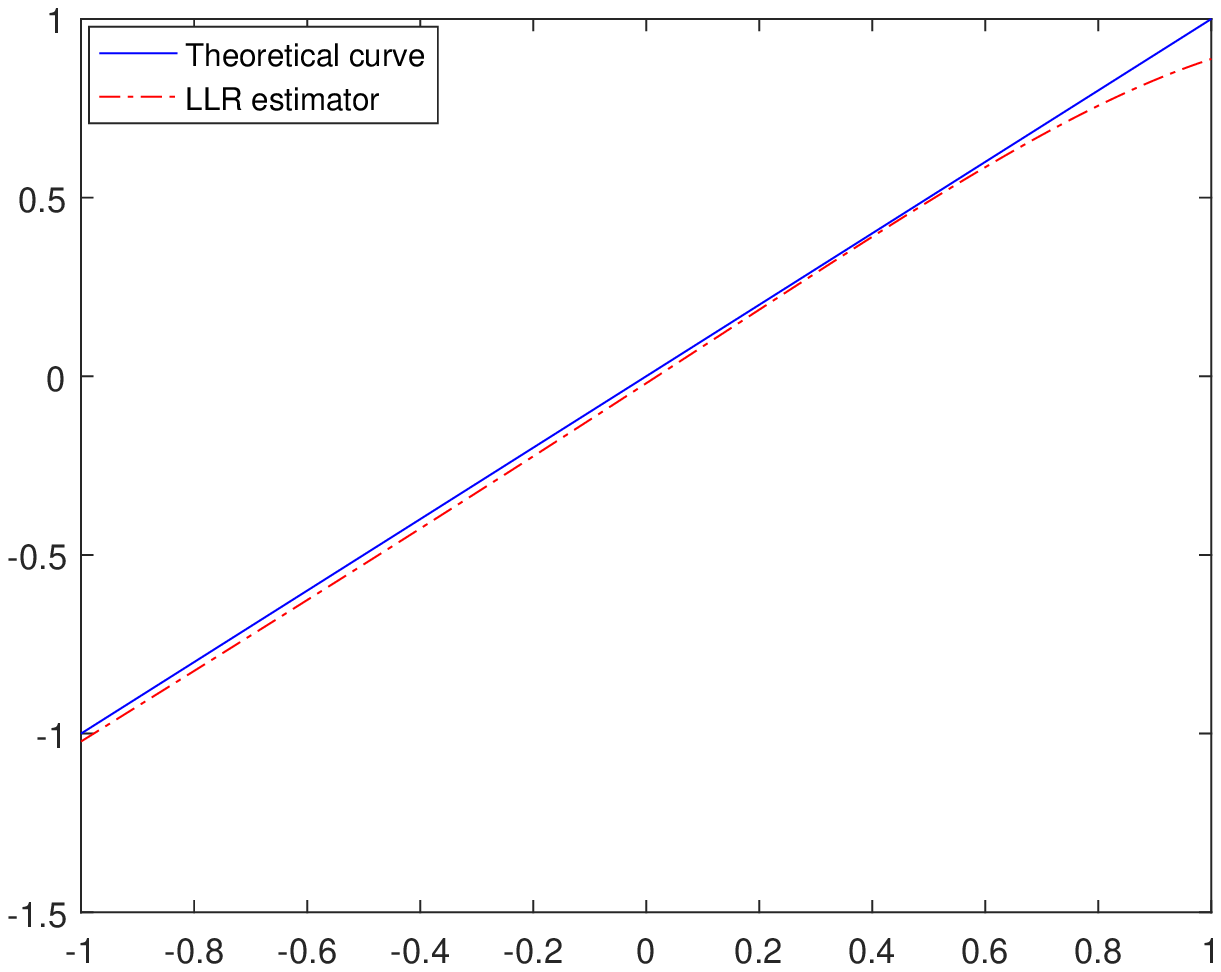}
	\end{minipage} \hfill
	\begin{minipage}[c]{.26\linewidth}
		\includegraphics[height=2in, width=2.3in]{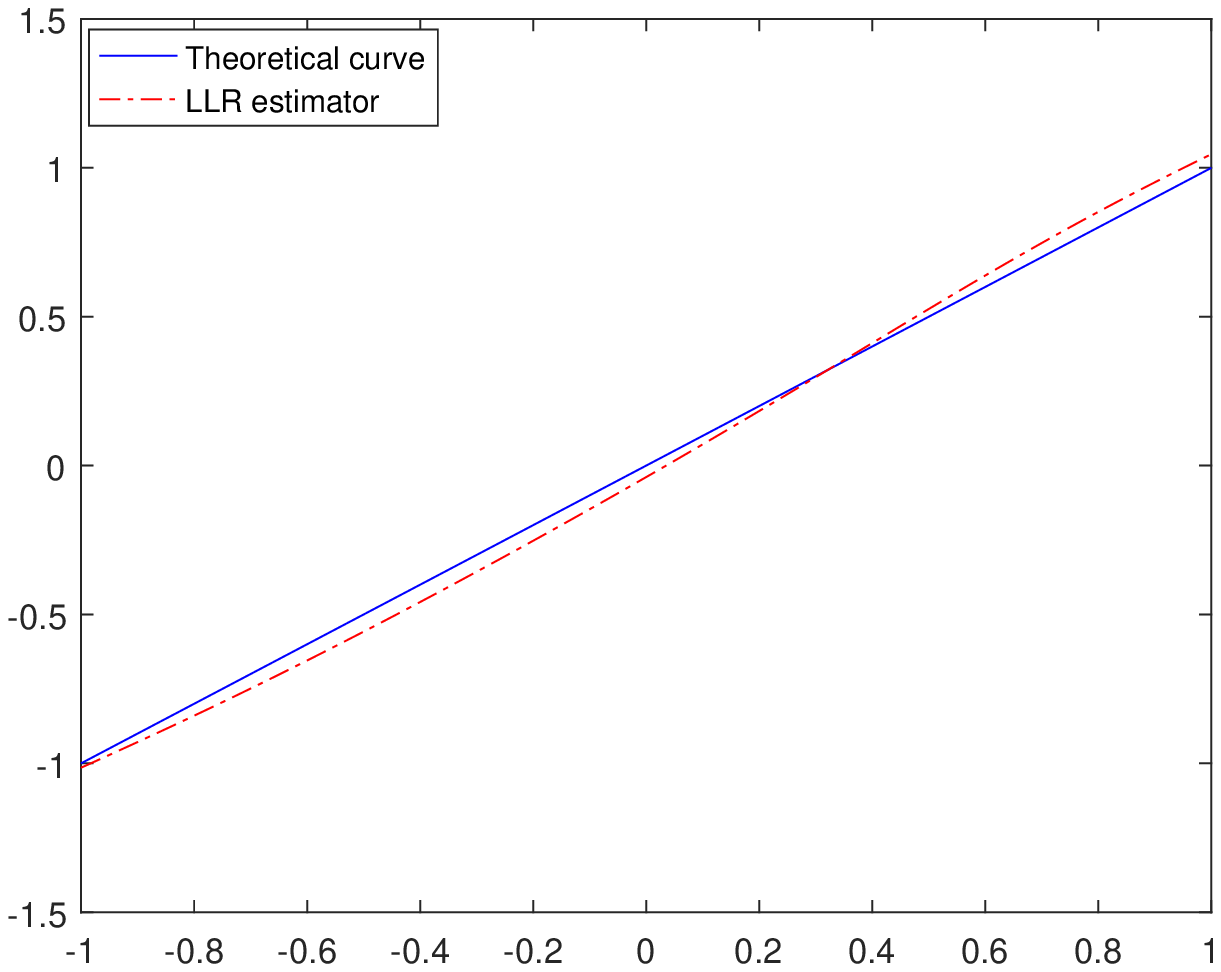}
	\end{minipage}\hfill\hfill
	\caption{\textcolor{blue}{$\quad \mu(\cdot)$}, \textcolor{red}{$\widehat{\mu}(\cdot)$} with $ n=300$ for C.P.$\approx 8, 25$ and $60\%$ respectively.}\label{figure2}
\end{figure}
%%%%%%%%%%%%%%%%
\begin{figure}[!h]
	\begin{minipage}[c]{.26\linewidth}
		\includegraphics[height=2in, width=2.3in]{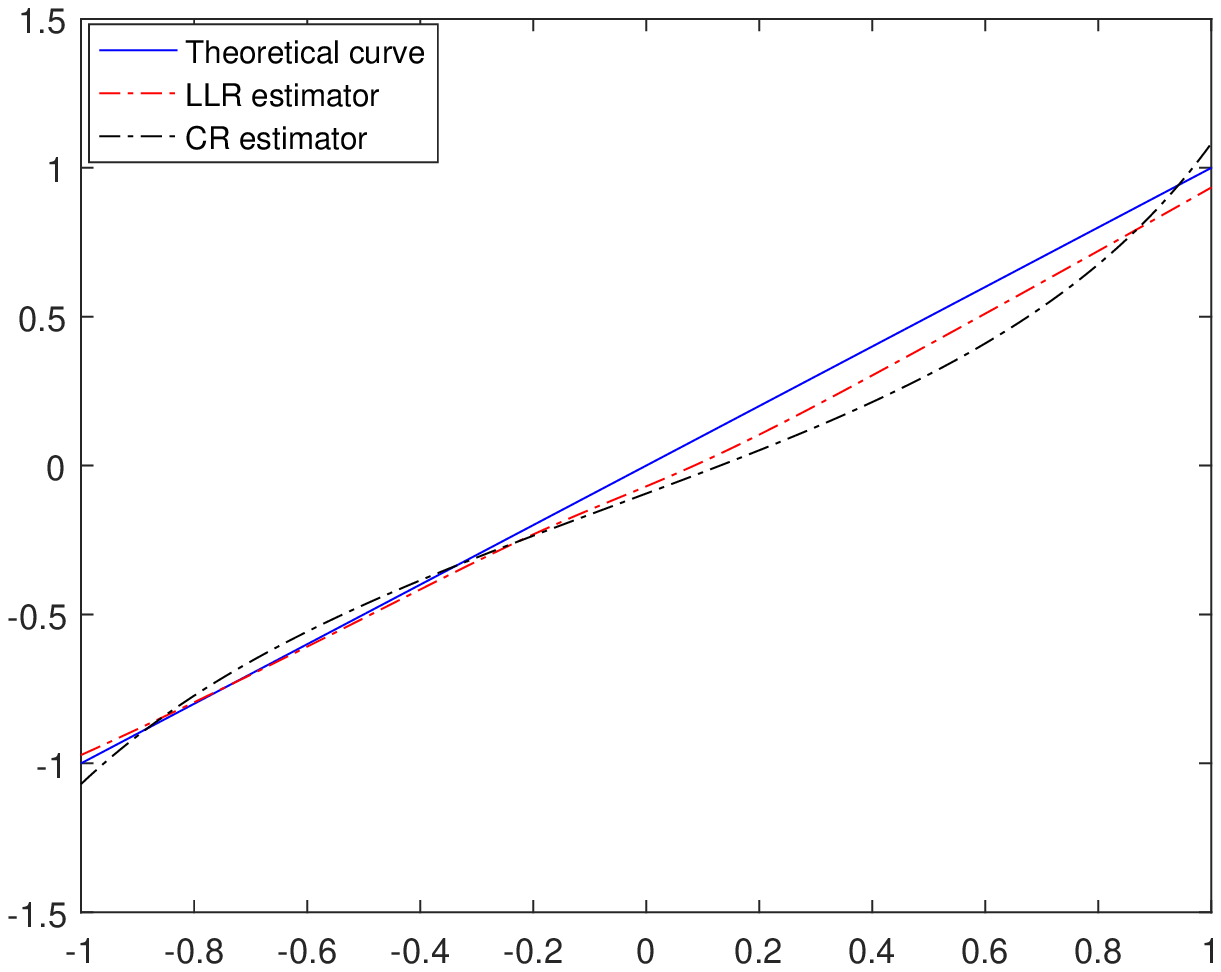}
	\end{minipage} \hfill
	\begin{minipage}[c]{.26\linewidth}
		\includegraphics[height=2in, width=2.3in]{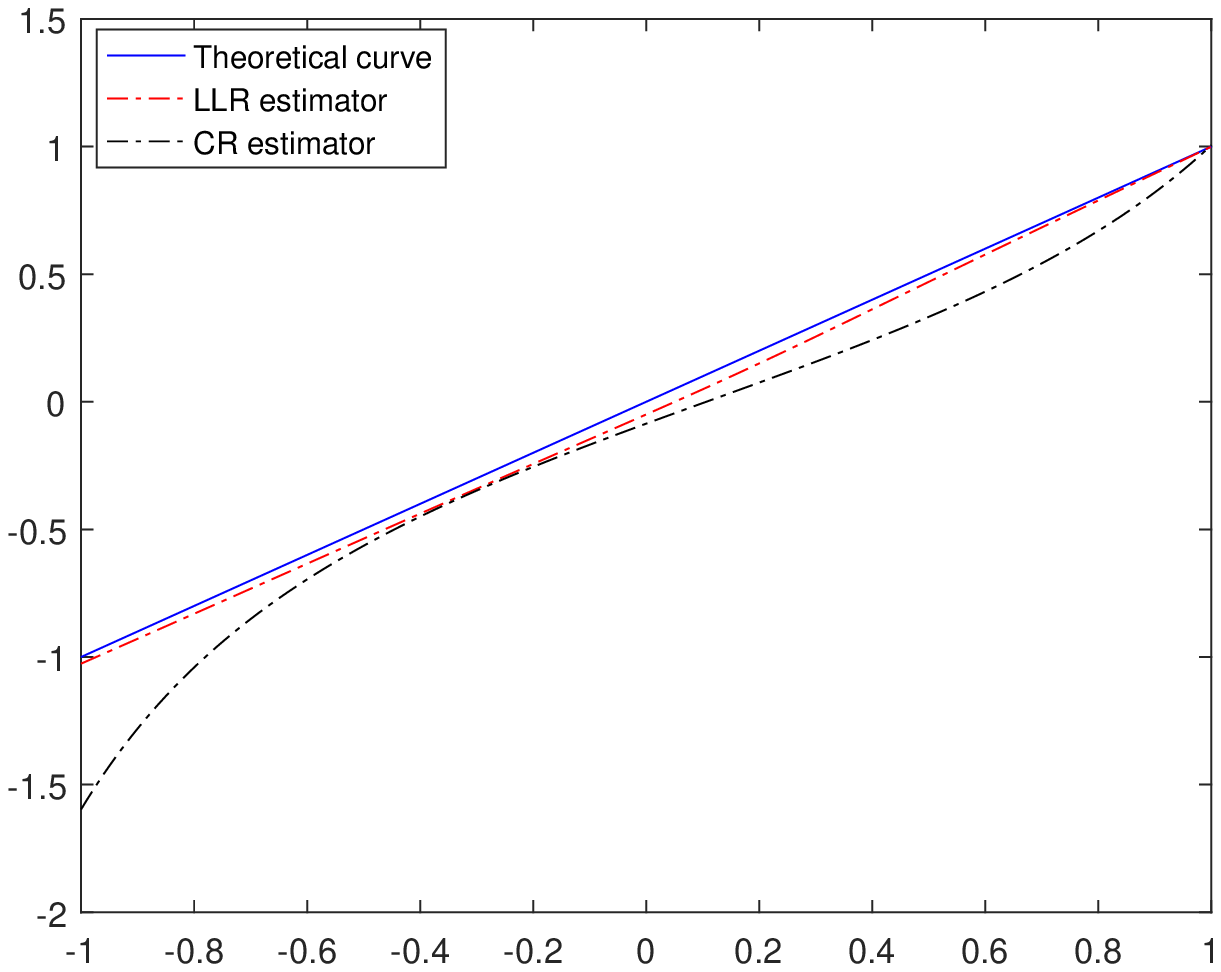}
	\end{minipage} \hfill
	\begin{minipage}[c]{.26\linewidth}
		\includegraphics[height=2in, width=2.3in]{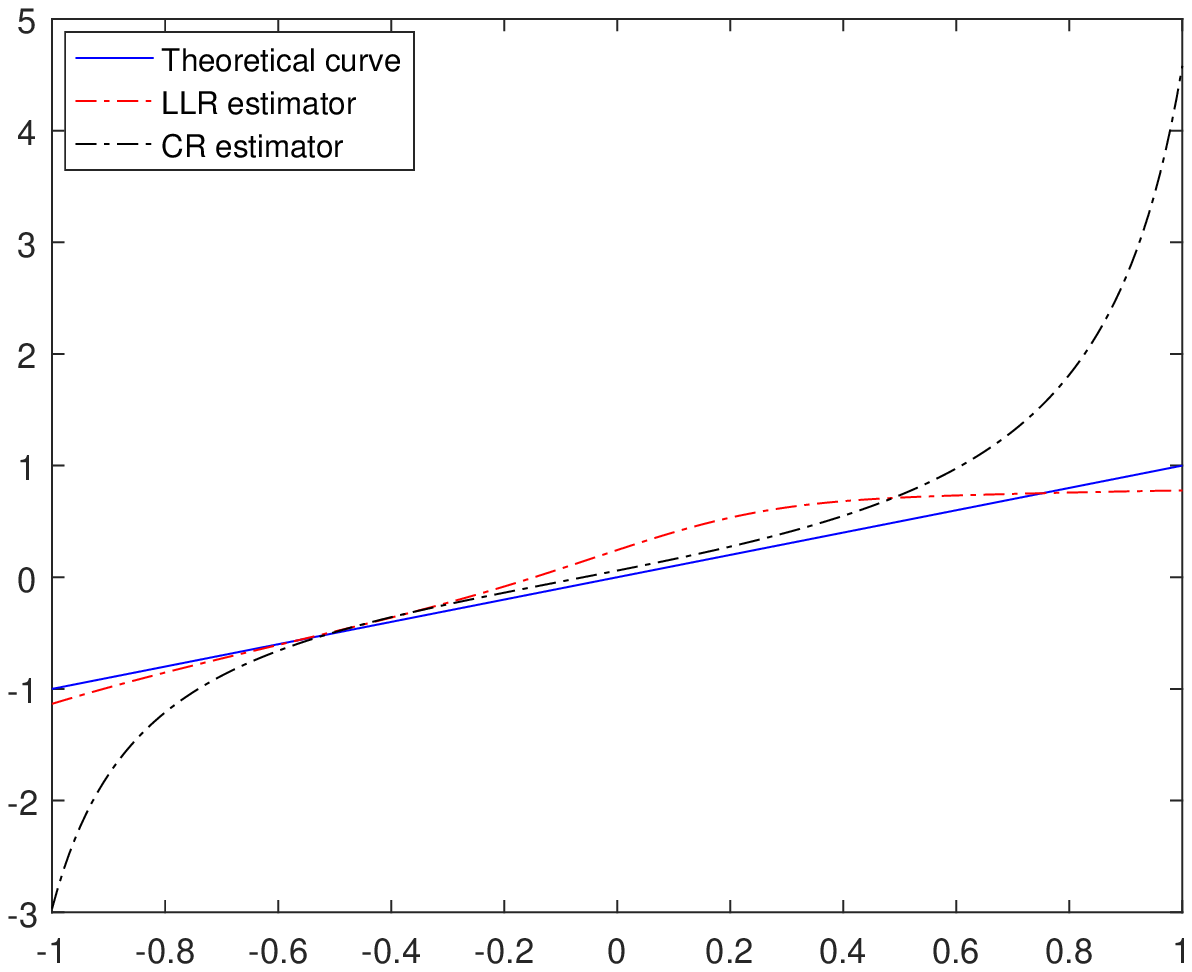}
	\end{minipage}\hfill\hfill
	\caption{\textcolor{blue}{$ \quad \mu(\cdot)$}, \textcolor{red}{$\widehat{\mu}(\cdot)$} and \textcolor{black}{$\mu_n(\cdot)$} with $ n=300$ for C.P.$\approx 10, 35$ and $65\%$ respectively.}\label{figure3}
\end{figure}
\begin{center}
	\begin{table}[t]%
		\centering
		\caption{Comparative table of MSE.\label{tab1}}%
		\begin{tabular*}{400pt}{@{\extracolsep\fill}lccc@{\extracolsep\fill}}%
			\hline
			\textbf{C.P.} & \textbf{n} & \textbf{CR} & \textbf{LLR} \\
			\hline
			& 100 & 0.0158 & 0.0011\\
			10 & 300 &0.0024 &0.0003\\
			&500 &2.48 $\times 10^{-4}$  &2.32 $\times 10^{-6}$\\
			\hline
			&100 & 0.0836 & 0.0025 \\
			30&300 &0.0473 & 0.0020\\
			&500 & 0.0108 & 8.10 $\times 10^{-4}$\\
			\hline
			&100 & 0.0611 & 0.1181\\
			50&300 & 0.2321 & 0.0228\\
			&500 & 0.0258 & 0.0064\\ 
			\hline
		\end{tabular*} 	
	\end{table}
\end{center}
In \hyperref[figure1]{Figure 1-3}, it can be seen that  for $(i)$ the LLR estimator performs better as increasing of the sample size $n$; $(ii)$ the estimator quality is affected by the C.P. but resists and keeps close to the theoretical curve; $(iii)$ the LLR and CR estimators are almost indistinguishable when the censorship rate is low. Notably, the CR estimator is sensitive to the effect of censorship, which is visible on the edges, unlike the LLR estimator, which resists the edge and remains stuck to the theoretical curve.
In the \hyperref[tab1]{Table 1}, we take different values of C.P. and report the mean squared error (MSE) of the LLR estimator and CR estimator. We can see that the LLR estimator performs better when the sample size increases and is only slightly affected by the percentage of observed data.
%%%%%%%---------- PROOFS ------------------------------------------%%%%%%%%%%%%%%%%%%%%%%%%%
\section{Proofs and Auxiliary results}\label{sect 5}
%%%%%%%%%%%%%%%%%%%%%%%%%%%%%%%%%%%%%%%%%%%%%%%%%%%
{\bf Proof of Proposition 1.3.} \label{proof of prop3}We retain all the notation from \hyperref[sect 2]{Section 2} and let denote by
	\begin{equation*}
	\widehat{S}_\ell(x)=\frac{1}{nh} \sum_{i=1}^n \widehat{Z}^\star_i (X_i-x)^\ell K_i, \quad \quad \widetilde{S}_\ell(x)=\frac{1}{nh} \sum_{i=1}^n {Z}^\star_i (X_i-x)^\ell K_i  \quad \text{for} \quad \ell=0,1,
	\end{equation*}
	and 
	\begin{equation*}
	   \widehat{T}_\ell(x)=\frac{1}{nh} \sum_{i=1}^n (X_i-x)^\ell K_i \quad \text{for} \quad \ell=0,1,2.
	\end{equation*} 
	Consider now the following decomposition: 
	%%%%%%%%%%%%%%%%%%%%%%%%%%%%%%%%%%%%%%%%%%%
	\begin{equation*}
	\begin{aligned}
	\left|\widehat{\mu}_1(x)-\widetilde{\mu}_1(x)\right|
	&=\left|\widehat{S}_0(x)\widehat{T}_2(x)-\widehat{S}_1(x)\widehat{T}_1(x)-\widetilde{S}_0(x)\widehat{T}_2(x)+\widetilde{S}_1(x)\widehat{T}_1(x)\right|\\
	&\leq \left|\widehat{T}_2(x)-\mathbb{E}[\widehat{T}_2(x)]\right|\times \left|\widehat{S}_0(x)-\widetilde{S}_0(x)\right|+\left|\mathbb{E}[\widehat{T}_2(x)]\left(\widehat{S}_0(x)-\widetilde{S}_0(x)\right)\right|\\
	&+\left|\widehat{T}_1(x)-\mathbb{E}[\widehat{T}_1(x)]\right|\times\left|\widehat{S}_1(x)-\widetilde{S}_1(x)\right|+\left|\mathbb{E}[\widehat{T}_1(x)]\left(\widehat{S}_1(x)-\widetilde{S}_1(x)\right) \right|,\\
	\end{aligned}
	\end{equation*}
	%%%%%%%%
	\noindent We then state and prove  \hyperref[lem1]{Lemma 1-3}  which are needed in \hyperref[prop3]{Proposition 1.3.} 
	%%%%%%%%%%%%%%%%%%%%% LEMMA 1%%%%%%%%%%%%%%%%%%%%%%%%%%%%%%%%%%%%%%
	\begin{lemma}\label{lem1}
		Under Assumptions \hyperref[A2]{A2} and \hyperref[A3]{A3}, we have for $\ell=0,1$
		\begin{equation*}
		\sup_{x \in \mathcal{C}} |\widehat{S}_{\ell}(x)-\widetilde{S}_{\ell}(x)|=\text{O}_{a.s.}\left(\sqrt{\frac{\log\log n }{n}}\right) \;\;\; \text{as}\;\;\; n \rightarrow \infty.
		\end{equation*}
	\end{lemma}
	%%%%%%%%%%%%%%%%%%%%%%%%%
	{\bf Proof of Lemma 1.}
		For $\ell=0,1,$ we have 
		\begin{equation*}
		\begin{aligned}
		\sup_{x \in \mathcal{C}}\left|\widehat{S}_\ell(x)-\widetilde{S}_\ell(x)\right|&= \sup_{x \in \mathcal{C}}\left| \frac{1}{n h} \sum_{j=1}^n \widehat{Z}_j^{\star} (X_j-x)^\ell K_j - \frac{1}{n h} \sum_{j=1}^n Z_j^{\star} (X_j-x)^\ell K_j \right|\\
		&\leq \sup_{x \in \mathcal{C}}\left| \frac{1}{n h} \sum_{j=1}^n Z_j (X_j-x)^\ell K_j \left( \frac{1}{\bar{G}_n(Z_j)}-\frac{1}{\bar{G}(Z_j)}\right)\right|\\
		&\leq \frac{1}{ \bar{G}^2(\tau_H)} \sup_{t \leq \tau_H} \left|\bar{G}_n(t)-\bar{G}(t)\right|\times \sup_{x \in \mathcal{C}}\left| \frac{1}{n h}\sum_{j=1}^n  Z_j (X_j-x)^\ell K_j\right|\\
		&=: \sup_{t \leq \tau_H} \mathcal{L}_1(t) \times \sup_{x \in \mathcal{C}} \left|\mathcal{L}_2(x)\right|.
		\end{aligned}
		\end{equation*}
		For $\mathcal{L}_1$, using Lemma 4.2. in \hyperref[deheuvels2000]{Deheuvels and Einmahl (2000)}, we get 
		\begin{equation}\label{L1}
		\sup_{t \leq \tau_H} \mathcal{L}_1(t) = \text{O}_{a.s.}\left(\sqrt{\frac{\log \log n}{n}}\right).
		\end{equation}
		For $\mathcal{L}_2$, under Assumptions \hyperref[A2]{A2} and  \hyperref[A3]{A3}, using the strong large law numbers,  change of variable and Taylor expansion around $x$, we have 
		\begin{eqnarray}\label{L2}
		\sup_{x \in \mathcal{C}}\left|\mathcal{L}_2(x)\right| &\leq& 	C \sup_{x \in \mathcal{C}} \left|\mathbb{E}\left[h^{-1} (X_1-x)^\ell K_1\right] \right|\nonumber\\
		&=& h^{\ell+1} \sup_{x \in \mathcal{C}} |f^{\prime}(x)| \int v^{\ell+1} K(v)dv.
		\end{eqnarray}
		Finally, combining the results in (\ref{L1}) and (\ref{L2}) concludes the proof of \hyperref[lem1]{Lemma 1}.
	%%%%%%%%%%%%%%%%%%%% LEMME 2 %%%%%%%%%%%%%%%%%%%%%%%%%
	\begin{lemma} \label{lem2}
		Under Assumption \hyperref[A1]{A1}, \hyperref[A2]{A2} and \hyperref[A3]{A3}, we have for $\ell=0,1,2$
		\begin{equation*}
		\sup_{x \in \mathcal{C}} \left|\widehat{T}_{\ell}(x)-\mathbb{E}[\widehat{T}_{\ell}(x)]\right|=\text{O}_{a.s.}\left(\sqrt{\frac{\log n }{n h }}\right) \;\;\; \text{as}\;\;\; n \rightarrow \infty.
		\end{equation*}
	\end{lemma}
	%%%%%%%%%%%%%%%%%%%%%%%%%%%%%% PROOF OF LEMME 2 %%%%%%%%%%%%
	{\bf Proof of Lemma 2.}
		Let $b_n=n^{-1/2\nu}$ for $\nu>0$ and cover the compact set $\mathcal{C}$ by $\displaystyle{ \cup_{i=1}^{d_n}}(x_i-b_n,x_i+b_n)$ with $d_n=\text{O}(n^{1/2\nu})$. 
	Let
		\begin{equation*}
		\mathcal{C}_n=\{x_i-b_n; \, x_i+b_n, \; 1 \leq i \leq d_n\},
		\end{equation*}
		the extremities of the latter subdivision, 
	Then
		\begin{equation} 
		\sup_{x \in \mathcal{C}} |\widehat{T}_{\ell}(x)-\mathbb{E}[\widehat{T}_{\ell}(x)]| \leq \max_{1\leq i \leq d_n} \max_{x \in \mathcal{C}_n} |\widehat{T}_{\ell}(x)-\mathbb{E}[\widehat{T}_{\ell}(x)]| + 2^{\nu} C b_n^{\nu}.
		\end{equation}
		Since  $b_n=n^{-1/2\nu}$ then  \[b_n^{\nu}=\text{O}\left(\sqrt{\frac{\log n}{nh}}\right).\] 
		Then for all $\varepsilon>0$, we have
		\begin{equation*}
		\mathbb{P}\left(\max_{x\in \mathcal{C}_n}|\widehat{T}_{\ell}(x)-\mathbb{E}[\widehat{T}_{\ell}(x)]|>\varepsilon \right) \leq \sum_{x \in \mathcal{C}_n} \mathbb{P} \left(|\widehat{T}_{\ell}(x)-\mathbb{E}[\widehat{T}_{\ell}(x)]|>\varepsilon\right).
		\end{equation*}
		Let us write for $\ell=0,1,2$ and $x \in \mathcal{C}_n$
		\begin{equation*}
		\begin{aligned}
		\widehat{T}_{\ell}(x)-\mathbb{E}[\widehat{T}_{\ell}(x)]&=\frac{1}{n} \sum_{i=1}^{n}\frac{(X_i-x)^{\ell} K_i - \mathbb{E}\left[(X_i-x)^{\ell} K_i\right]}{h}\\
		&=:\frac{1}{n} \sum_{i=1}^{n} \mathcal{A}_{\ell,i}(x).
		\end{aligned}
		\end{equation*}
		Since an almost complete property holds almost surely, we apply \hyperref[coro]{Corollary 1}. For that, we focus on the absolute moments of $\mathcal{A}_{\ell,i}(x)$
		\begin{equation*}
		\begin{aligned}
		\mathbb{E}|\mathcal{A}_{\ell,i}(x)|^m&= h^{-m} \mathbb{E} \left| \sum_{k=0}^{m} c_{k,m} \left((X_i-x)^{\ell} K_i\right)^k \mathbb{E}\left[(X_i-x)^{\ell} K_i\right]^{m-k}\right|\\
		&\leq h^{-m} \sum_{k=0}^{m} c_{k,m}  \left| \mathbb{E}\left[\left((X_1-x)^{\ell} K_1\right)^k\right] \mathbb{E}\left[(X_1-x)^{\ell} K_1\right]^{m-k}\right|.\\
		\end{aligned}
		\end{equation*}
		On the one hand, using conditional expectation property, \hyperref[A2]{A2}  and \hyperref[A3]{A3}, we have
		\begin{equation*}
		\begin{aligned}
		\mathbb{E}\left[\left((X_1-x)^{\ell} K_1\right)^k\right]&=h^{\ell k+1}\int v^{\ell k} K^k(v)f(x+vh)dv
		\end{aligned}
		\end{equation*}
		and 
		\begin{equation*}
		\begin{aligned}
		\mathbb{E}\left[(X_1-x)^{\ell} K_1\right]^{m-k}	&=\left(h^{\ell+1}\int v^{\ell} K(v)f(x+vh)dv\right)^{m-k}\\
		\end{aligned}
		\end{equation*}
		then, for $\ell=0,1,2$ and $\forall m \geq 2$
		\begin{equation*}
		\begin{aligned}
		\mathbb{E}|\mathcal{A}_{\ell,1}(x)|^m &= \text{O}(h^{m\ell-k+1})
		&= \text{O}(\max_{1 \leq k \leq m}h^{-k+1})
		&= \text{O}(h^{-m+1}).
		\end{aligned}
		\end{equation*}
		We can now apply \hyperref[coro]{Corollary 2}. Choosing $a^2=h^{-1}$ we get
		\begin{equation*}
		\begin{aligned}
		\mathbb{P}\left( \left| \widehat{T}_{\ell}(x)-\mathbb{E}[\widehat{T}_{\ell}(x)] \right|>\varepsilon \right)&=\mathbb{P}\left( \left|\sum_{i=1}^n \mathcal{A}_{\ell,i}(x) \right|>\varepsilon n\right)\\
		&\leq 2 \exp \left( -\frac{\varepsilon^2 nh}{2(1+\varepsilon)}\right).
		\end{aligned}
		\end{equation*}
		Hence, for $\varepsilon=\varepsilon_0 \left(\frac{\log n}{nh}\right)^{1/2}$ and $n$ large enough, we get
		\begin{equation*}
		\mathbb{P}\left( \left| \widehat{T}_{\ell}(x)-\mathbb{E}[\widehat{T}_{\ell}(x)] \right|>\varepsilon\right) \leq 2 \exp \left(-\frac{\varepsilon_0^2}{4} \log n\right)=2 n^{-\frac{\varepsilon_0^2}{4}}.
		\end{equation*}
		It follows that 
		\begin{equation*}
		\sum_{x \in \mathcal{C}_n} \mathbb{P}\left(|\widehat{T}_{\ell}(x)-\mathbb{E}[\widehat{T}_{\ell}(x)]|>\varepsilon\right) \leq 4 n ^{-\frac{\varepsilon_0^2}{4}+\frac{1}{2\nu}}.
		\end{equation*}
		Finally, an appropriate choice of $\varepsilon_0$ yields an upper bound of order $n^{-3/2}$ which by Borel-Cantelli's lemma completes the proof of \hyperref[lem2]{Lemma 2}.
	%%%%%%%%%%%%% LEMMA 3 %%%%%%%%%%%%%%%%%%%%%%%%%%%%
	\begin{lemma}\label{lem3}
		Under Assumptions  \hyperref[A1]{A1}, \hyperref[A2]{A2} and \hyperref[A3]{A3}, we have for $\ell=0,1,2$
		\begin{equation*}
		\mathbb{E}[\widehat{T}_{\ell}(x)]=\text{O}(h^{\ell}) .\;\;\;
		\end{equation*}
	\end{lemma}
	%%%%%%%%%%%%%%%%%%%% Proof of Lemma 3 %%%%%%%%%%%%%%%%%%%%%%%
	{\bf Proof of Lemma 3.}
		For $\ell=0,1,2$, using a change of variable and Taylor expansion for $\xi \in [x, x+hv]$, we have 
		\begin{equation*}
		\begin{aligned}
		\mathbb{E}[\widehat{T}_{\ell}(x)]	&=  h^{\ell} \int v^{\ell} K(v) f(x+hv)dt\\
		&= h^{\ell} f(x) \int v^{\ell} K(v) dv + h^{\ell+1} \int v^{\ell+1} K(v) f^{\prime}(\xi) dv\\
		\end{aligned}
		\end{equation*}
		under Assumptions \hyperref[A1]{A1}, \hyperref[A2i]{A2} and \hyperref[A3]{A3} we get the result.
	%%%%%%%%%%%%%%%%%%%%%%%%
	\noindent Then, combining the results in \hyperref[lem1]{Lemma 1} with \hyperref[lem2]{Lemma 2} and \hyperref[lem1]{Lemma 1} with \hyperref[lem3]{Lemma 3} we get the result of \hyperref[prop1]{Proposition 1}.
%%%%%%%%%%%%%%%%%%%%%%%%%%%%% Proof of Proposition 4 %%%%%%%%%%%%%%%%
{\bf Proof of Proposition 1.4.} By similar reasoning than the  \hyperref[proof of prop3]{Proof of Proposition 1.3}, we remark that:
	\begin{equation*}
	\begin{aligned}
	\mathcal{B}_{4}(x)
	&=\left\{\widetilde{S}_0(x)\widehat{T}_2(x)-\mathbb{E}\left[\widetilde{S}_0(x)\widehat{T}_2(x)\right]\right\}-\left\{\widetilde{S}_1(x)\widehat{T}_1(x)-\mathbb{E}\left[\widetilde{S}_1(x)\widehat{T}_1(x)\right] \right\}\\
	&=: \mathcal{B}_{4,1}(x)-\mathcal{B}_{4,2}(x).
	\end{aligned}
	\end{equation*}
	\noindent On the one hand 
	\begin{eqnarray}\label{decom1}
	\mathcal{B}_{4,1}(x)&=&\left(\widetilde{S}_0(x)-\mathbb{E}\left[\widetilde{S}_0(x)\right]\right)\left(\widehat{T}_2(x)-\mathbb{E}\left[\widehat{T}_2(x)\right]\right)+\left(\widehat{T}_2(x)-\mathbb{E}\left[\widehat{T}_2(x)\right]\right)\mathbb{E}\left[\widetilde{S}_0(x)\right]\nonumber\\
	&+&\left(\widetilde{S}_0(x)-\mathbb{E}\left[\widetilde{S}_0(x)\right]\right)\mathbb{E}\left[\widehat{T}_2(x)\right]+\mathbb{E}\left[\widetilde{S}_0(x)\right]\mathbb{E}\left[\widehat{T}_2(x)\right]-\mathbb{E}\left[\widetilde{S}_0(x)\widehat{T}_2(x)\right].
	\end{eqnarray}
	\noindent On the other hand 
	\begin{eqnarray}\label{decom2}
	\mathcal{B}_{4,2}(x)&=&\left(\widetilde{S}_1(x)-\mathbb{E}\left[\widetilde{S}_1(x)\right]\right)\left(\widehat{T}_1(x)-\mathbb{E}\left[\widehat{T}_1(x)\right]\right)+\left(\widehat{T}_1(x)-\mathbb{E}\left[\widehat{T}_1(x)\right]\right)\mathbb{E}\left[\widetilde{S}_1(x)\right]\nonumber\\
	&+&\left(\widetilde{S}_1(x)-\mathbb{E}\left[\widetilde{S}_1(x)\right]\right)\mathbb{E}\left[\widehat{T}_1(x)\right]
	+\mathbb{E}\left[\widetilde{S}_1(x)\right]\mathbb{E}\left[\widehat{T}_1(x)\right]-\mathbb{E}\left[\widetilde{S}_1(x)\widehat{T}_1(x)\right].
	\end{eqnarray}
	It remains to study each term in the decompositions (\ref{decom1}) and (\ref{decom2}). 
	%All the terms  $\widehat{T}_\ell(x)-\mathbb{E}\left[\widehat{T}_{\ell}(x)\right]$ and $\mathbb{E}\left[\widehat{T}_{\ell}(x)\right]$ for $\ell=1,2$ were considered in \hyperref[lem2]{Lemma 2} and \hyperref[lem3]{Lemma 3} respectively. For the others terms, 
	For that, let us consider the following Lemmas.
	%%%%%%%%%%%%%%%%%%%%%%%%%%%%%%%%%
	%%%%%%%%%%%%%%%%%%%%%%%%%%%%%%%%%%%%% LEMMA 4 %%%%%%%%%%%%%%%%%%%
	\begin{lemma}\label{lem4}
		Under Assumption \hyperref[A1]{A1}, \hyperref[A2]{A2}, \hyperref[A4]{A4} and \hyperref[A5]{A5}, we have for $\ell=0,1,$
		\begin{equation*}
		\sup_{x \in \mathcal{C}}\left|\widetilde{S}_{\ell}(x)-\mathbb{E}\left[\widetilde{S}_{\ell}(x)\right]\right|=\text{O}_{a.s.}\left(\sqrt{\frac{\log n}{n\,h}}\right) \quad \quad \quad \text{as} \quad n \rightarrow \infty.
		\end{equation*}
	\end{lemma}
	%%%%%%%%%%%%%%%%%%%%%%%%%%%%%%%%% PROOF OF lEMMA 4 %%%%%%%%%%%%%%%%%%%%
	{\bf Proof of Lemma 4.}
		The proof is very similar to that of \hyperref[lem2]{Lemma 2}. The same notations of \hyperref[lem2]{Lemma 2} are used. 
		\[\sup_{x \in \mathcal{C}} |\widetilde{S}_{\ell}(x)-\mathbb{E}[\widetilde{S}_{\ell}(x)]| \leq \max_{1\leq j \leq d_n} \max_{x \in \mathcal{C}_n} |\widetilde{S}_{\ell}(x)-\mathbb{E}[\widetilde{S}_{\ell}(x)]| + 2^{\nu} C b_n^{\nu}.\]
		Observe that 
		\begin{equation*}
		\mathbb{P}\left(\max_{x\in \mathcal{C}_n}|\widetilde{S}_{\ell}(x)-\mathbb{E}[\widetilde{S}_{\ell}(x)]|>\varepsilon \right) \leq \sum_{x \in \mathcal{C}_n} \mathbb{P} \left(|\widetilde{S}_{\ell}(x)-\mathbb{E}[\widetilde{S}_{\ell}(x)]|>\varepsilon\right).
		\end{equation*}	
		Let us write for $\ell=0,1$ and $x \in \mathcal{C}_n$
		\begin{equation*}
		\begin{aligned}
		\widetilde{S}_{\ell}(x)-\mathbb{E}[\widetilde{S}_{\ell}(x)]&=\frac{1}{n} \sum_{j=1}^{n}\frac{Z^{\star}_j(X_j-x)^{\ell} K_j - \mathbb{E}\left[Z^{\star}_j(X_j-x)^{\ell} K_j\right]}{h}\\
		&=:\frac{1}{n} \sum_{j=1}^{n} \mathcal{R}_{\ell,j}(x).
		\end{aligned}
		\end{equation*}
		In order to apply \hyperref[coro]{Corollary 2}, we focus on the absolute moments of $\mathcal{R}_{\ell,j}(x)$ for $\ell=0,1$
		\begin{equation*}
		\begin{aligned}
		\mathbb{E}|\mathcal{R}_{\ell,j}(x)|^m &\leq h^{-m} \sum_{k=0}^{m} c_{k,m}  \left| \mathbb{E}\left[\left(Z^{\star}_1(X_1-x)^{\ell} K_1\right)^k\right] \mathbb{E}\left[Z^{\star}_1(X_1-x)^{\ell} K_1\right]^{m-k}\right|.\\
		\end{aligned}
		\end{equation*}
		On the one hand, using  conditional expectation property, under \hyperref[A5]{A5}  for $m \geq k$, we get
		\begin{eqnarray}\label{synt_cal}
		\mathbb{E}\left[\left(Z^{\star}_1(X_1-x)^{\ell} K_1\right)^k\right]&=&\mathbb{E}\left[(X_1-x)^{\ell k} K_1^k \mathbb{E}[Z^{\star,k}_1|X_1]\right]\nonumber\\ 
		&=& \int (u-x)^{\ell k } K^{k}\left(\frac{u-x}{h}\right) \mathbb{E}[Z^{\star,k}_1|X_1=u]f_X(u)du
		\end{eqnarray}
		with
		\begin{eqnarray}\label{synt_cal2}
		\mathbb{E}[Z^{\star,k}_1|X_1=u]&=&\mathbb{E}\left[\frac{Z_1^k}{\bar{G}^{k-1}(Y_1)}\big|X_1=u\right]\nonumber\\
		&\leq& \frac{1}{\bar{G}^{k-1}(\tau_H)} \int z^k f_{Z|X}(z|u)dz.
		\end{eqnarray}
		We replace (\ref{synt_cal2}) in (\ref{synt_cal}), under \hyperref[A5]{A5}, we get 
		\begin{equation*}
		\begin{aligned}
		\mathbb{E}\left[\left(Z^{\star}_1(X_1-x)^{\ell} K_1\right)^k\right]
		&=\mathbb{E}\left[(X_1-x)^{\ell k} K_1^k \mathbb{E}[Z^{\star,k}_1|X_1]\right]\\ 
		&=\int (u-x)^{\ell k } K^{k}\left(\frac{u-x}{h}\right) \mathbb{E}[Z^{\star,k}_1|X_1=u]f_X(u)du\\
		&\leq \frac{1}{\bar{G}^{k-1}(\tau_H)} \int(u-x)^{\ell k} K^k\left(\frac{u-x}{h}\right) \int z^k f_{Z|X}(z|u)f_X(u)dzdu \\
		&= \frac{1}{\bar{G}^{k-1}(\tau_H)} \int(u-x)^{\ell k} K^k\left(\frac{u-x}{h}\right) \upsilon_k(u) du \\
		&= \frac{h^{\ell k+1}}{\bar{G}^{k-1}(\tau_H)} \int s^{\ell k} K^k(s) \upsilon_k(x+hs) ds. \\
		\end{aligned}
		\end{equation*}
		\noindent    On the other hand, under (\ref{synt}) and analogously to the previous development, we get 
		\begin{equation*}
		\begin{aligned}
		\mathbb{E}\left[Z^{\star}_1(X_1-x)^{\ell} K_1\right]^{m-k}
		&=\mathbb{E}\left[(X_1-x)^{\ell} K_1\mathbb{E}[Z^{\star}_1|X_1]\right]^{m-k}\\
		&=\mathbb{E}\left[(X_1-x)^{\ell} K_1\mathbb{E}\left[Z_1|X_1\right]\right]^{m-k}\\
		&=\left(\int (u-x)^{\ell} K\left(\frac{u-x}{h}\right)\mu(u)f(u)du\right)^{m-k}\\
		&=\left(h^{\ell+1}\int v^{\ell} K(v)S_0(x+vh)dv\right)^{m-k}.\\
		\end{aligned}
		\end{equation*}
		By \hyperref[A2]{A2} and \hyperref[A4]{A4}, for $\ell=0,1$ and $\forall \;m \geq 2$
		\begin{equation*}
		\begin{aligned}
		\mathbb{E}|\mathcal{R}_{\ell,1}(x)|^m &\leq \text{O}(h^{m\ell-k+1})
		&= \text{O}(\max_{1 \leq k \leq m}h^{-k+1})
		&= \text{O}(h^{-m+1}).
		\end{aligned}
		\end{equation*}
		Now, we can apply \hyperref[coro]{Corollary 2}. By choosing $a^2=h^{-1}$, $\varepsilon=\varepsilon_0 \left(\frac{\log n}{nh}\right)^{1/2}$ and for $n$ large enough, we get
		\begin{equation*}
		\begin{aligned}
		\mathbb{P}\left( \left| \widetilde{S}_{\ell}(x)-\mathbb{E}[\widetilde{S}_{\ell}(x)] \right|>\varepsilon\right)&=\mathbb{P}\left( \left|\sum_{j=1}^n \mathcal{R}_{\ell,j}(x) \right|>n\varepsilon\right)\\
		& \leq 2 \exp \left(-\frac{\varepsilon_0^2}{4} \log n\right)=2 n^{-\frac{\varepsilon_0^2}{4}}.
		\end{aligned}
		\end{equation*}
		It follows that 
		\begin{equation*}
		\sum_{x \in \mathcal{C}_n} \mathbb{P}\left(|\widetilde{S}_{\ell}(x)-\mathbb{E}[\widetilde{S}_{\ell}(x)]|>\varepsilon\right) \leq 4 n ^{-\frac{\varepsilon_0^2}{4}+\frac{1}{2\nu}}.
		\end{equation*}
		Finally, an appropriate choice of $\varepsilon_0$ yields to an upper bound of order $n^{-3/2}$ which by Borel-Cantelli's lemma completes the proof of \hyperref[lem4]{Lemma 4}.
	%%%%%%%%%%%%%%%%%%%%%%%%%%%%%%%%%%%%%%%%% LEMME 5 %%%%%%%%%%%%%%%%%%%
	\begin{lemma}\label{lem5}
		Under Assumptions \hyperref[A1]{A1}, \hyperref[A2]{A2} and \hyperref[A4]{A4}, we have for $\ell=0,1,$
		\begin{equation*}
		\mathbb{E}[\widetilde{S}_{\ell}(x)]=\text{O}(h^{\ell})\quad\quad \quad \text{as} \quad n \rightarrow \infty.
		\end{equation*}
	\end{lemma}
	%%%%%%%%%%%%%%%%%%%%%%%%%%%%%%%%%%%%% PROOF OF LEMMA 5 %%%%%%%%%%%%%%%%
	{\bf Proof of Lemma 5.}
		\noindent For $\ell=0,1$, using a change of variable, Taylor expansion for $\xi \in ]x,x+vh[$ and under \hyperref[A1]{A1}, \hyperref[A2]{A2} and \hyperref[A4]{A4}, we have
		\begin{equation*}
		\begin{aligned}
		\mathbb{E}\left[\widetilde{S}_{\ell}(x)\right]&=\int (hv)^{\ell}K(v)S_0(x+hv)dv\\
		&=h^{\ell} S_0(x) \int v^{\ell} K(v)dv+h^{\ell+1} \int v^{\ell+1} K(v) S_0^{\prime}(\xi)dv.
		\end{aligned}
		\end{equation*}
	\noindent Now, it remains to study the quantity $\mathbb{E}\left[\widetilde{S}_0(x)\right]\mathbb{E}\left[\widehat{T}_2(x)\right]-\mathbb{E}\left[\widetilde{S}_0(x)\widehat{T}_2(x)\right]$. To do that, let consider the following Lemma.
	%%%%%%%%%%%%%%%%%%%%%%%%%%%%%% lEMMA 6 %%%%%%%%%%%%%%%%%%%%%%%%%%%%%%%%%%%%
	\begin{lemma}\label{lem6}
		Under Assumptions \hyperref[A1]{A1-A4}, we have 
		\begin{equation*}
		\text{Cov}(\widetilde{S}_0(x),\widehat{T}_2(x))=\text{o}\left(\sqrt{\frac{\log n}{nh}}\right) \quad \quad \quad \text{as} \quad n \rightarrow \infty.
		\end{equation*}
	\end{lemma}
	%%%%%%%%%%%%%%%% PROOF OF LEMMA 6 %%%%%%%%%%%%%%%%%%%%%%%%%%%%%%%%%%%%%%%%
	{\bf Proof of Lemma 6.} By a change of variable, Taylor expansion and under \hyperref[A1]{A1-A4}, we have 
		\begin{equation*}
		\begin{aligned}
		\text{Cov}(\widetilde{S}_0(x),\widehat{T}_2(x))
		&= \frac{1}{(n h)^2} \sum_{j=1}^n \sum_{i=1}^n \left[\mathbb{E}\left(Z_j^{\star}(X_i-x)^2 K_iK_j\right)-\mathbb{E}\left(Z_j^{\star} K_j\right)\mathbb{E}\left((X_i-x)^2 K_i\right) \right]\\
		&=\left(\frac{n(n-1)-n^2}{(nh)^2}\right)\mathbb{E}\left((X_1-x)^2 K_1\right)\mathbb{E}\left(Z_1^{\star} K_1\right)+\frac{n}{(nh)^2}\mathbb{E}\left(Z_1^{\star} (X_1-x)^2 K_1^2\right)\\
		&=O\left(\frac{h}{n}\right)
		\end{aligned}
		\end{equation*}
		which is negligible with respect to $\sqrt{\frac{\log n}{nh}}$.
	%%%%%%%%%%%%%%%%%%%%%%%%%%%%% LEMMA 7 %%%%%%%%%%%%%%%%%%%%%%%%%%%%%%%%%%%%%%%%%%%%
	\begin{lemma}\label{lem7}
		Under Assumptions \hyperref[A1]{A1-A4}, we have 
		\begin{equation*}
		\text{Cov}(\widetilde{S}_1(x),\widehat{T}_1(x))=\text{o}\left(\sqrt{\frac{\log n}{nh}}\right)\quad \quad \quad \text{as} \quad n \rightarrow \infty.
		\end{equation*}
	\end{lemma}
	%%%%%%%%%%%%%%%%%%%%%%%%% PROOF OF LEMMA 7 %%%%%%%%%%%%%%%%%%%%%%%%%%%%%%%%%%%
	{\bf Proof of Lemma 7.}
		By a change of variable, Taylor expansion and under\hyperref[A1]{A1-A4}, we have 
		\begin{equation*}
		\begin{aligned}
		\text{Cov}(\widetilde{S}_1(x),\widehat{T}_1(x))
		&= \frac{1}{(n h)^2} \sum_{j=1}^n \sum_{i=1}^n \left[\mathbb{E}\left(Z_j^{\star}(X_i-x)(X_j-x) K_iK_j\right)-\mathbb{E}(Z_j^{\star} K_j)\mathbb{E}\left((X_i-x) K_i\right) \right]\\
		&=\left(\frac{n(n-1)-n^2}{(nh)^2}\right)\mathbb{E}\left((X_1-x) K_1\right)\mathbb{E}\left(Z_1^{\star} (X_1-x) K_1\right)+\frac{n}{(nh)^2}\mathbb{E}\left(Z_1^{\star} (X_1-x)^2 K_1^2\right)\\
		&=O\left(\frac{h^2}{n}\right)
		\end{aligned}
		\end{equation*}
		which is negligible with respect to $\sqrt{\frac{\log n}{nh}}$.
	\noindent Then combining the results in \hyperref[lem2]{Lemma 2} and \hyperref[lem4]{Lemma 4} with \hyperref[lem2]{Lemma 2} and \hyperref[lem5]{Lemma 5} with \hyperref[lem3]{Lemma 3} and \hyperref[lem4]{Lemma 4} in addition to \hyperref[lem6]{Lemma 6} and \hyperref[lem7]{Lemma 7} concludes the proof of the \hyperref[prop4]{Proposition 1.4}.\\
%%%%%%%%%%%%%%%%%%%%%%% Proof of Proposition 2 %%%%%%%%%%%%%%%%%%%
{\bf Proof of Proposition 1.2.}
	This proof is similar to that of \hyperref[prop4]{Proposition 1.4}. We use the same notation used in \hyperref[lem2]{Lemma 2} and the following decomposition 
	\begin{equation*}
	\begin{aligned}
	\mathcal{B}_{2}(x)&=\widehat{T}_0(x)\widehat{T}_2(x)-\widehat{T}_1^2(x)-\mathbb{E}[\widehat{T}_0(x)\widehat{T}_2(x)-\widehat{T}_1^2(x)]\\
	&=\left\{\widehat{T}_0(x)\widehat{T}_2(x)-\mathbb{E}[\widehat{T}_0(x)\widehat{T}_2(x)]\right\}-\left\{\widehat{T}_1^2(x)-\mathbb{E}[\widehat{T}_1^2(x)]\right\}\\
	&=:\mathcal{B}_{2,1}(x)-\mathcal{B}_{2,2}(x) .
	\end{aligned}
	\end{equation*}	
	On the one hand 
	\begin{eqnarray}\label{B21}
	\mathcal{B}_{2,1}(x)&=&(\widehat{T}_0(x)-\mathbb{E}[\widehat{T}_0(x)])(\widehat{T}_2(x)-\mathbb{E}[\widehat{T}_2(x)])+\mathbb{E}[\widehat{T}_0(x)](\widehat{T}_2(x)-\mathbb{E}[\widehat{T}_2(x)])\nonumber\\
	&+&\mathbb{E}[\widehat{T}_2(x)](\widehat{T}_0(x)-\mathbb{E}[\widehat{T}_0(x)])+\mathbb{E}[\widehat{T}_2(x)]\mathbb{E}[\widehat{T}_0(x)]-\mathbb{E}[\widehat{T}_0(x)\widehat{T}_2(x)].
	\end{eqnarray}
	On the other hand 
	\begin{equation}\label{B22}
	\mathcal{B}_{2,2}(x)=\widehat{T}_1^2(x)-\mathbb{E}[\widehat{T}_1^2(x)]=\text{Var}[\widehat{T}_1(x)].
	\end{equation}
	it remains to study each term in (\ref{B21}) and (\ref{B22}). The terms $\mathbb{E}[\widehat{T}_{\ell}(x)]$ and $\widehat{T}_{\ell}(x)-\mathbb{E}[\widehat{T}_{\ell}(x)]$ for $\ell=0,1,2$ were considered in \hyperref[lem3]{Lemma 3} and \hyperref[lem2]{Lemma 2} respectively. For the others terms, we consider the following Lemmas.
	%%%%%%%%%%%%%%%%%%%%%%%%%% LEMMA 8 %%%%%%%%%%%%%%%%%%%%%%%%%%%%%%%
	\begin{lemma}\label{lem8}
		Under Assumptions \hyperref[A1]{A1}, \hyperref[A2]{A2} and \hyperref[A3]{A3}, for $\ell=0,1,2$, we have  
		\begin{equation*}
		\text{Var}[\widehat{T}_{\ell}(x)]=\text{o}\left(\sqrt{\frac{\log n}{n h}}\right)\quad\quad \quad \text{as} \quad n \rightarrow \infty.
		\end{equation*}
	\end{lemma}
	%%%%%%%%%%%%%%%%%%%%% PROOF OF LEMMA 8 %%%%%%%%%%%%%%%%%%%%%%%%%%%%%% 
	{\bf Proof of Lemma 8.}
		For $\ell=0,1,2,$ 	with a change of variable and Taylor expansion with $\xi \in ]x,x+hs[$, we have 
		\begin{equation*}
		\begin{aligned}
		\mathbb{E}[\widehat{T}_{\ell}^2(x)]
		&=\frac{1}{nh^2}\int (u-x)^{2\ell}K^2\left(\frac{u-x}{h}\right) f(u)du\\
		&=\frac{1}{nh}\int (sh)^{2\ell}K^2(s) f(x+hs)ds\\ 
		&=\frac{1}{nh}\left\{ h^{2\ell} f(x) \int s^{2\ell} K^2(s) ds+ h^{2\ell+1} \int s^{2\ell+1} K^2(s)  f^{\prime}(\xi)ds\right\}.
		\end{aligned}
		\end{equation*}	
		Under \hyperref[A1]{A1}, \hyperref[A2]{A2} and \hyperref[A3]{A3}  we get the result. Furthermore the result is $\text{o}\left(\sqrt{\frac{\log n}{nh}}\right)$.
	%%%%%%%%%%%%%%%%%%%%%%%%%%%%%%%%% LEMMA 9 %%%%%%%%%%%%%%%%%%%%%%%%%%%%%%%%%%%
	\begin{lemma}\label{lem9}
		Under Assumptions \hyperref[A1]{A1-A3}, we have 
		\begin{equation*}
		\text{Cov}(\widehat{T}_0(x),\widehat{T}_2(x))=\text{o}\left(\sqrt{\frac{\log n}{nh}}\right) \quad\quad \quad \text{as} \quad n \rightarrow \infty.
		\end{equation*}
	\end{lemma}
	%%%%%%%%%%%%%%%%%%%%%%%% PROOF OF LEMMA 9 %%%%%%%%%%%%%%%%%%%%%%%%%%% 
	{\bf Proof of Lemma 9.}
		By a change of variable and Taylor expansion we have 
		\begin{equation*}
		\begin{aligned}
		\text{Cov}(\widehat{T}_0(x),\widehat{T}_2(x))
		&= \frac{1}{(n h)^2} \sum_{j=1}^n \sum_{i=1}^n \left[\mathbb{E}\left((X_i-x)^{2} K_iK_j \right)-\mathbb{E}\left( K_j\right)\mathbb{E}\left((X_i-x)^{2} K_i\right) \right]\\
		&= \frac{1}{n^2h^2} \left\{(n(n-1)-n^2)\mathbb{E}\left[ K_1\right]\mathbb{E}\left[(X_1-x)^{2} K_1\right]+\mathbb{E}\left[(X_1-x)^{2} K^2_1 \right]\right\}\\
		&= O\left(\frac{h^2}{n}\right).
		\end{aligned}
		\end{equation*}
		Under Assumptions \hyperref[A1]{A1-A3}, we get the result.
	%%%%%%%%%%%%%%%%%%%%%%%%%%%%%%%%%%%%%%%%%%%%%%%%%%%%%%
	\noindent Then, combining the results in \hyperref[lem2]{Lemma 2} and \hyperref[lem2]{Lemma 2} with \hyperref[lem3]{Lemma 3} in addition to the results in \hyperref[lem8]{Lemma 8} and \hyperref[lem9]{Lemma 9} conclude the proof of the \hyperref[prop2]{Proposition 1.2}.\\
	%%%%%%%%%%%%%%%%%%%%%%%%%%%%%% Proof of corrollary 1 %%%%%%%%%%%%%%%%%%%%
	\noindent{\bf Proof of corollary 1.} There exists $\Gamma>0$, such that for all $ x \in \mathcal{C}, \;\E[\widehat{\mu}_0(x)] \geq \Gamma$.
	Therefore $\displaystyle \inf_{x \in \mathcal{C}}\widehat{\mu}_0(x) \leq \frac{\Gamma}{2} $ implies that there exists $x \in \mathcal{C}$ such that $|\E[\widehat{\mu}_0(x)]-\widehat{\mu}_0(x)|\geq \frac{\Gamma}{2}$ which gives 
	\begin{equation*}
	\sup_{x \in \mathcal{C}}\Big|\E[\widehat{\mu}_0(x)]-\widehat{\mu}_0(x)\Big|\geq \frac{\Gamma}{2}.
	\end{equation*}
	Thus, the result of \hyperref[prop2]{Proposition 1.2.} allows to write that for $\frac{\Gamma}{2}=\Gamma^{\prime}$:
	\begin{equation*}
	\sum_{n} \mathbb{P} \left( \inf_{x \in \mathcal{C}} \widehat{\mu}_0(x) \leq \Gamma^{\prime} \right) \leq \sum_{n} \mathbb{P} \left( \sup_{x \in \mathcal{C}} \Big|\E[\widehat{\mu}_0(x)]-\widehat{\mu}_0(x)\Big| \leq \Gamma^{\prime} \right)<\infty.
	\end{equation*}
%%%%%%%%%%%%%%%%%%%%%%%%%%%%%%%%%%%% Proposition 1 %%%%%%%%%%%%%%%%%%
{\bf Proof of Proposition 1.1.} Let consider 
	\begin{equation*}
	\begin{aligned}
	|\mathcal{B}_1(x)|=\left| \frac{\mathbb{E}\left[\widetilde{\mu}_1(x)\right]}{\mathbb{E}\left[\widehat{\mu}_0(x)\right]}-\mu(x)\right|
	&=\left| \frac{\mathbb{E}\left[\widetilde{\mu}_1(x)\right]-\mu(x)\mathbb{E}\left[\widehat{\mu}_0(x)\right]}{\mathbb{E}\left[\widehat{\mu}_0(x)\right]}\right|\\
	&=\left|\frac{h^{-2}\left\{\mathbb{E}[w_{1,2}(x)Z^{\star}_2]-\mu(x)\mathbb{E}\left[w_{1,2}(x)\right]\right\}}{h^{-2}\mathbb{E}\left[w_{1,2}(x)\right]}\right|\\
	&=\left|\frac{\mathbb{E}\left[w_{1,2}(x)\left\{\mathbb{E}[Z^{\star}_2|X_2]-\mu(x)\right\}\right]}{\mathbb{E}[w_{1,2}(x)]}\right|\\
	&=\left|\frac{\mathbb{E}\left[w_{1,2}(x)\left(\mu(X_2)-\mu(x)\right)\right]}{\mathbb{E}\left[w_{1,2}(x)\right]}\right|\\
	&=\left|\mu(X_2)-\mu(x)\right|.\\
	\end{aligned}
	\end{equation*}
	Thus, under \hyperref[A6]{A6}, we have 
	\[\sup_{x \in \mathcal{C}} |\mathcal{B}_1(x)|\leq C|X_2-x|^{\nu} \leq Ch^{\nu}.\]
\noindent Finally, by summing the results in \hyperref[prop1]{Proposition 1.1}-\hyperref[prop1]{Proposition 1.4}, we get the proof of \hyperref[theo1]{Theorem 1}.
%------------------ Corollary ---------------------%%%%%%%%%%%%%%%
\begin{coro}{(A.8. p. 234 in \hyperref[Ferraty2006]{Ferraty and Vieu (2006)})} \label{coro}. Let $U_i$ be a sequence of independent r.v. with zero mean. If $\forall \; m \geq 2$, $\exists \; C_m>0$, $\mathbb{E}[|U_1^m|] \leq C_m a^{2(m-1)}$, we have 
	\begin{equation*}
	\forall \varepsilon>0, \;\;\; \mathbb{P}\left( \left|\sum_{i=1}^{n} U_i\right|> n\varepsilon\right) \leq 2 \exp \left\{-\frac{\varepsilon^2 n}{2a^2(1+\varepsilon)}\right\}.
	\end{equation*}
\end{coro}
%%%%%%%%%%%   ----------- Concluding remarks ----------------------%%%%%%%%%%%%%
\section*{Concluding remarks} 
In this note, we study a local linear  regression  function estimator and show that this method has advantages with respect to  the classical kernel estimator. On simulated data we show that the LL method is more efficient than the classical kernel method. On the one hand, it mitigates edge effects and on the other hand, it remains efficient when the censoring rate increases substantially. We point out that the method proposed in \hyperref[cai2003]{Cai (2003)} with two weights where the first is standard kernel for smoothing and the second which is the \hyperref[kaplan1958]{Kaplan and Meier (1958)} estimator. The resulting estimator is interesting, however the author does not use a weighting in the denominator with the survival law of the censoring random variable. 
%The latter reduces the effect of censorship at beginning and at the end of the study in particular if the rate of censoring is very high. 
Recall that \hyperref[ElGhouch2008]{El Ghouch and Van Keilegom (2008)} estimated the regression function by applying polynomial local linear regression techniques using Beran's estimator. Their conditions need to have a result about conditional law on the censored random variable that in our case we do not use it.  Furthermore their uniform result is given only  in probability. We point out and to the best of our knowledge, the type of our result  has never been obtained.


\begin{thebibliography}{}
	\bibitem{}\label{beran1981}
	Beran, R. (1981). Nonparametric regression with randomly censored survival data, Tech. Report, University of California,
	\bibitem{}\label{cai2003}
	Cai, Z. (2003).  Weighted local linear approach to censored nonparametric regression. \emph{Recent Advances and Trends in Nonparametric Statist. Michael G. Akritas and Dimitris N. Politis (Editors)}.
	\bibitem{}\label{carbonez1995}
	Carbonez A., Gyorfi L., Van Der Meulen E.C. (1995). Partitioning estimates of a regression function under random censoring. {\it Statist. and Decisions.} {\bf 76}, 1335--1344.
	\bibitem{}\label{dabrowska1987}
	Dabrowska, D. (1987). Nonparametric regression with censored survival data. {\it Scand. J. Statist.,} {\bf 14}, 181--197. 
	\bibitem{}\label{deheuvels2000}
	Deheuvels, P., Einmahl, J.H.J. (2000). Functional limit laws for the increments of Kaplan-Meier product-limit processes and applications. {\it The Annals of Probability.} {\bf 28,} 1301--1335.
	\bibitem{}\label{ElGhouch2008}  
	El-Ghouch, M., Van Keilegom, I. (2008). Nonparametric regression with dependent censored data. \textit{Scand. J. Statist.} \textbf{35(2)}, 228--247.
	\bibitem{}\label{ElGhouch2009}  
	El-Ghouch, M., Van Keilegom, I. (2009). Local linear quantile regression with dependent censored data. {\it Statist. Sinica} \textbf{19}, 1621--1640.
	\bibitem{}\label{fan1992}
	Fan, J. (1992). Design-adaptative Nonparametric Regression. \emph{J. of the American Statist. Association.} {\bf 87,} 998--1004.
	\bibitem{}\label{fan1994}
	Fan, J. and Gijbels, I.(1994). Censored regression: local linear approximation and their applications. {\it J. of the Ameri. Statist. Assoc.} {\bf 89(426)}. 560--570.
	\bibitem{}\label{fan1996}
	Fan, J., Gijbels, I. (1996). {\it Local Polynomial Modeling and Its Applications.} {\bf  66}.  Chapman \& Hall/CRC.
	\bibitem{}\label{fan2003}
	Fan, J., Yao, Q. (2003). \emph{Nonlinear time series: nonparametric and parametric methods}. Springer, New York.
	\bibitem{}\label{Ferraty2006}
	Ferraty, F., Vieu, P. (2006). {\it  Nonparametric Functional Data Analysis: Theory and Practice.}  \emph{Springer}. New York.
	\bibitem{}\label{Guessoum2008}
	Guessoum, Z., Ould Sa\"\i d E. (2008).
	On nonparametric estimation of the regression function under random censorship model. \emph{Statist. and Decisions,} {\bf 26},1001--1020.
	\bibitem{}\label{kim1998}
	Haesook T., Kim, Young K. Truong. (1998).
	Nonparametric regression estimates with censored data: local linear smoothers and their applications. \emph{Biometrics.} {\bf 54.} 1434--1444.
	\bibitem{}\label{kaplan1958} 
	Kaplan, E.L., Meier, P. (1958). Nonparametric estimation from incomplete observations. \emph{ J. Amer. Stat. Assoc.}, {\bf  53}, 458--481.
	\bibitem{}\label{kohler2002}
	K\"ohler, M., M\'{a}th\`e, K., Pint\"er, M. (2002). Prediction from randomly right censored data. \emph{J. Multivar. Anal.,} {\bf  80}, 73--100.
	\bibitem{}\label{Nadaraya1964}
	Nadaraya, E. A. (1964). On estimating regression. \emph{Theor.  Probab. Appl.} {\bf 9}, 141--142.
	\bibitem{}\label{Watson1964}
	Watson, G.S. (1964). Smooth regression analysis. \emph{ Sankhy\`a}, A, {\bf 26}, 359--372.
\end{thebibliography}
\end{document}